\documentclass{amsart}

\usepackage{latexsym}
\usepackage{amscd}
\usepackage{amssymb}
\usepackage{amsmath}
\usepackage{xy}
\xyoption{all}
\usepackage{amsthm}
\usepackage{enumerate}
\usepackage{hyperref}
\usepackage{graphics}

\newtheorem{thm}{Theorem}[section]
\newtheorem{lem}[thm]{Lemma}
\newtheorem{defn}[thm]{Definition}
\newtheorem{prop}[thm]{Proposition}
\newtheorem{cor}[thm]{Corollary}
\newtheorem{eg}[thm]{Example}
\newtheorem{rmk}[thm]{Remark}

\numberwithin{equation}{section}

\newcommand{\leftexp}[2]{{}^{#1}{#2}}

\def\Z{\mathbf{Z}}

\def\P{\mathbf{P}}
\def\del{\partial}

\def\ra{\rightarrow}

\def\del{\partial}

\def\coh{\operatorname{coh}}
\def\cohc{\operatorname{coh,c}}
\def\Coh{\operatorname{Coh}}

\def\perf{\operatorname{perf}}

\def\id{\operatorname{id}}
\def\Id{\operatorname{Id}}

\def\Hom{\operatorname{Hom}}

\def\vect{\operatorname{vect}}
\def\hocolim{\operatorname{hocolim}}

\begin{document}

\author{Matthew Robert Ballard}
\address{Department of Mathematics, University of Pennsylvania, Philadelphia, PA, USA}
\email{ballardm@math.upenn.edu}

\begin{abstract}
 We extend Orlov's result on representability of equivalences to schemes projective over a field. We also investigate the quasi-projective case.
\end{abstract}

\title[Equivalences of derived categories]{Equivalences of derived categories of sheaves on quasi-projective schemes}
\maketitle

\section{Introduction}

Given two projective schemes, $X$ and $Y$, and an exact functor, $F: D^b_{\coh}(X) \ra D^b_{\coh}(Y)$, one can wonder if, on the product $X \times Y$, there exists a bounded complex coherent sheaves, $E$, whose associated integral transform is isomorphic to $F$. We can try to compare $D^b_{\coh}(X \times Y)$ with the category of exact functors between $D^b_{\coh}(X)$ and $D^b_{\coh}(Y)$ by studying the functor that takes a complex on the product to the associated integral transform. One's first hope might be that this functor is an equivalence. A simple look at the case of the projective line shatters this hope. The space of natural endomorphisms of the identity functor of $D^b_{\coh}(\P^1_k)$ is infinite-dimensional (assuming $k$ is infinite) while $D^b_{\coh}(\P^1_k \times \P^1_k)$ has finite-dimensional morphisms. (The author thanks Chris Brav for pointing this out.) Moreover, any morpshim $\phi: E \ra F[2]$ between coherent sheaves $E$ and $F$ on $\P^1_k \times \P^1_k$ induces the zero natural transformation. Thus, the morphism sets of the two categories are very different. However, being stubborn as we are, this does not extinquish the flame of our hope; it only focuses it.

The next best scenario would be that any exact functor, $F: D^b_{\coh}(X) \ra D^b_{\coh}(Y)$, is isomorphic to integral transform and the kernel of the integral transform is unique up to quasi-isomorphism. This hope is more difficult to stamp out. Indeed, there are no counterexamples and there is some supporting evidence. The evidence comes from Orlov in the form of the following theorem, see \cite{Orl97}:
\begin{thm}
 Let $X$ and $Y$ be smooth projective varieties and let $F: D^b_{\coh}(X) \ra D^b_{\coh}(Y)$ be a full and faithful exact functor possessing either a left or right adjoint. Then, there exists a bounded complex of coherent sheaves, $E$, on $X \times Y$, whose associated integral transform, $\Phi_E$, restricted to $D^b_{\coh}(X)$, is isomorphic to $F$.
\end{thm}
In particular, an important case covered by Orlov's result is where $F$ is an equivalence. 

This paper is an attempt to relax the conditions of smoothness and projectivity in the hypotheses of Orlov's theorem, and, consequently, provide more evidence in favor of a bijection between exact functor on derived categories and objects on the product up to quasi-isomorphism. It contains a few results on projective schemes that are quite similar to Orlov's result. For instance, we have the following theorem:

\begin{thm}
 Let $X$ and $Y$ be projective schemes over a field $k$. If $F: D_{\perf}(X) \ra D_{\perf}(Y)$ is a full and faithful functor with left and right adjoints, then $F$ is isomorphic to the restriction of $\Phi_E$ to $D_{\perf}(X)$ for an $E \in D^b_{\coh}(X \times Y)$.
\label{thm:1}
\end{thm}

And a corollary of it:

\begin{cor}
 Let $X$ and $Y$ be projective schemes over a field $k$. If $F: D^b_{\coh}(X) \ra D^b_{\coh}(Y)$ is an exact equivalence, then there exists an $E \in D^b_{\coh}(X \times Y)$ and a natural isomorphism $F \cong \Phi_E|_{D^b_{\coh}(X)}$.
\label{cor:1}
\end{cor}

These results are natural extensions Orlov's original result. However, in the intervening decade, little knowledge about the singular case has arisen. Even when other results, such as sufficient conditions for an integral transform to be an equivalence, are extended to more singular schemes, the question of whether derived-equivalent schemes are related by a Fourier-Mukai transform is side-stepped, see \cite{RMS07} and especially section $4$ of \cite{RMS06}. The results of this paper allow one to assign a kernel to any equivalence. One can then use geometric reasoning to study the kernel, a process which proves fruitful in the case of smooth and projective varieties.

If one wants to relax the projectivity assumption, we have the following:

\begin{thm}
 Let $X$ and $Y$ be quasi-projective schemes over a field $k$ with either $X$ or $Y$ possessing an ample line bundle, large tensor powers of which have trivial higher cohomology. If $D(X)$ and $D(Y)$ are equivalent, then there is an object of $D(X\times Y)$ whose associated integral transform is an equivalence.
\label{thm:2}
\end{thm}

Note one does not necessarily have isomorphism of the two functors appearing in theorem \ref{thm:2}. Indeed, the simplest case to check that such an isomorphism exists would be when both $X$ and $Y$ are affine. We jump wholly into the realm of commutative algebra and, surprisingly, we find no answers. The state of the knowledge remains essentially unchanged since Rickard's paper on derived Morita theory, \cite{Ric89}.

If one wishes to extend the results in this paper, there are a few obvious cases: stacks and twisted derived categories. However, from work Dugger and Shipley, see \cite{DS07}, if we push too far, equivalences are no longer guaranteed to arise from ``bimodules.'' It would be quite interesting to paint a line separating the Dugger-Shipley realm from the happy land presented here.

If one is willing to enrich the derived category by remembering more structure, i.e. the structure of a differential graded category or a stable $\infty$-category lying above the triangulated structure, one can show that all functors, preserving this extra structure, are integral transforms, see \cite{BLL04,Toe07,BFN08}. 

One can combine the results on dg-categories and dg-functors with work of Lunts and Orlov, \cite{LO09}. Independently, Lunts and Orlov prove a slightly stronger result involving projective scheme; one can remove the assumption of the existence of adjoints from \ref{thm:1}. They also prove a similar for bounded derived categories of coherent sheaves on projective schemes. These are applications of a central new idea: lifting structure from the triangulated category to a dg-enhancement. Whereas the methods of this paper seem to be difficult to adapt to a general exact functor, one can hope that Lunts and Orlov's results might be more amenable.

Here is a outline of the paper. In section \ref{sec:prelim}, we recall the results of \cite{Bal09a} which serve as the main new ingredient. In section \ref{sec:intrans}, we study integral transforms generally. We focus on the interplay between the existence of adjoints and the preservation of certain subcategories of the unbounded derived category of quasi-coherent sheaves. Section \ref{sec:tot} recalls Orlov's notion of convolution of a complex over a triangulated category. We mention how one can extend Orlov's ideas to totalize unbounded complexes via homotopy colimits of convolutions. In section \ref{sec:ample}, we recall another of Orlov's useful definitions: ample sequences in derived categories. We focus on ample sequences consisting of perfect objects. We talk about resolutions of the diagonal in section \ref{sec:resolvediag}. In the final two sections, we use the ideas and results of the previous sections to prove new results. In section \ref{sec:derivedMorita}, we discuss quasi-projective schemes, in particular proving the theorem \ref{thm:2}, In section \ref{sec:repfunct}, we focus on the projective case and prove theorem \ref{thm:1} amongst other results.

This work was a portion of the author's thesis at the University of Washington. The author would like to thank his advisor, Charles Doran, for his patience, energy, and dedication. While preparing this paper, the author was supported by NSF Research Training Group Grant, DMS 0636606.

\section{Preliminaries}
\label{sec:prelim}

Some notional preliminaries: given a category, $\mathcal{C}$, the morphism set from an object, $A$, to an object, $B$, is denoted as $[A,B]$. If a category is endowed with shift functor, the shift is denoted by $[1]$.

Given a scheme $X$, the category of perfect complexes, $D_{\perf}(X)$, is the full subcategory of the unbounded derived category of quasi-coherent sheaves on $X$, $D(X)$, consisting of complexes locally quasi-isomorphic to bounded complexes of finite rank locally-free sheaves. If $X$ is quasi-compact and separated, the objects of $D_{\perf}(X)$ are the compact objects of $D(X)$, meaning the natural map,
\begin{displaymath}
 \bigoplus_{i \in I} [A,B_i] \ra [A,\coprod_{i \in I} B_i],
\end{displaymath}
is an isomorphism for any perfect $A$ and any collection $B_i$, \cite{Nee92}. Any nonzero object of $D(X)$ admits a nonzero morphism from a compact object. Because of this, we call $D(X)$ compactly-generated. Brown's theorem on representability of cohomological functors on the category of spectra can be extended to compactly-generated triangulated categories, see \cite{Nee96}. It provides a useful tool for studying $D(X)$.

If we restrict to the situation where $X$ is quasi-projective over a field $k$, a complex is perfect if and only if it is globally quasi-isomorphic to a bounded complex of finite rank locally-free sheaves. While the notion of a perfection of an object is manifestly geometric, it is often not as useful as the more natural notion of compactness. Consequently, the identification of $D_{\perf}(X)$ as the subcategory of compact objects is helpful. Many other subcategories of $D(X)$ are defined geometrically; one can ask for a more intrinsic characterization of these subcategories. One such category is the bounded derived category of coherent sheaves with proper support, $D_{\cohc}^b(X)$. The desired characterization comes from viewing any object, $B \in D^b_{\cohc}(X)$, as a functor, $[-,B]: D_{\perf}(X)^{\circ} \ra \vect_k$, that takes triangles to long exact sequences and satsifies the following finiteness condition:
\begin{displaymath}
 \sum_{j \in \Z} \dim_k [A[j],B] < \infty
\end{displaymath}
for any perfect $A$. Any functor, $\phi: D_{\perf}(X)^{\circ} \ra \vect_k$, taking triangles to long exact sequences and satisfying the finiteness condition is called a locally-finite cohomological functor.
\begin{thm}
 For any quasi-projective scheme over a field, the bounded derived category of coherent sheaves with proper support is equivalent to the category of locally-finite cohomological functors. Moreover, the equivalence is the functor given above.
\end{thm}
For a proof, see \cite{Bal09a}. This result succeeds in providing the requested intrinsic characterization of $D^b_{\cohc}(X)$. It is in some sense dual to the category of compact objects. The duality is rather strong. It is easy to check that any object $A$ satisfying
\begin{displaymath}
 \sum_{j \in \Z} \dim_k [A[j],B] < \infty
\end{displaymath}
for all $B$ in $D^b_{\cohc}(X)$ must be perfect. Moreover, an impressive representability result of Rouquier, see \cite{Rou03}, immediately implies the following proposition.
\begin{prop}
 Let $X$ be a projective scheme over a perfect field $k$. The category of locally-finite homological functors on $D^b_{\coh}(X)$ is equivalent to $D_{\perf}(X)$ via the Yoneda embedding.
\label{prop:doubledual}
\end{prop}

This duality at the level of categories, including morphisms, allows us access to a wider range of tools than just an observation concerning objects would allow. For instance, given two quasi-projective schemes, $X$ and $Y$, and a functor, $F: D_{\perf}(X) \ra D_{\perf}(Y)$, there exists a unique functor, $F^{\vee}: D^b_{\cohc}(Y) \ra D^b_{\cohc}(X)$, with natural isomorphisms
\begin{displaymath}
 [FA,B] \cong [A,F^{\vee}B]
\end{displaymath}
for any $A$ in $D_{\perf}(X)$ and $B$ in $D^b_{\cohc}(Y)$. The proof of the existence is a simple application of the previous theorem. Uniqueness is clear. We call $F^{\vee}$ a right pseudo-adjoint to $F$. We also call $F$ left pseudo-adjoint to $F^{\vee}$. For a given functor $G: D^b_{\cohc}(Y) \ra D^b_{\cohc}(X)$, the existence of a left pseudo-adjoint, $\leftexp{\vee}{G}$, is guaranteed in the case of the proposition \ref{prop:doubledual} above. There are other cases where the existence of a left pseudo-adjoint is known. For instance, if $G: D^b_{\coh}(Y) \ra D^b_{\coh}(X)$ already possesses a left adjoint, $\leftexp{\vee}{G}$, $\leftexp{\vee}{G}$ must take perfect objects to perfect objects. Specializing further, if $G$ is an equivalence, then its inverse is its left pseudo-adjoint.

With this knowledge fresh in our memory, we begin the investigation in earnest by studying integral transforms homologically.

\section{Integral transforms}
\label{sec:intrans}

Let $X$ and $Y$ be quasi-compact, separated schemes and let $f: X \ra Y$ be a morphism.

\begin{defn}
 An object, $E$, from $D(X)$ is \textbf{$f$-perfect} if $f_*(E \otimes -):D(X) \ra D(Y)$ sends perfect objects to perfect objects.
\end{defn}

Recall that $f_*$ possesses a right adjoint $f^!$, \cite{Nee96}.

\begin{lem}
 $E$ is $f$-perfect if and only if $\mathcal{H}om(E,f^!-): D(Y) \ra D(X)$ commutes with coproducts.
\end{lem}

\proof $f_*(E \otimes -)$ is left adjoint to $\mathcal{H}om(E,f^!-)$ so the result follows from the next lemma. \qed

\begin{lem}
 Let $\mathcal{T}$ be a compactly-generated triangulated category and $\mathcal{S}$ a triangulated category. Let $F: \mathcal{T} \ra \mathcal{S}$ be a functor which commutes with coproducts and let $G: \mathcal{S} \ra \mathcal{T}$ be the right adjoint. $G$ commutes with coproducts if and only if $F$ takes a generating set of compact objects to compact objects.
\end{lem}

\proof Assume that $G$ commutes with coproducts. Let $X$ be a compact object of $\mathcal{T}$. Then,
\begin{displaymath}
 \bigoplus [F(X),Y_j] \cong \bigoplus [X,G(Y_j)] \cong [X,G(\coprod Y_j)] \cong [F(X),\coprod Y_j].
\end{displaymath}
The resulting map agrees with the natural map $\bigoplus [F(X),Y_j] \ra [F(X),\coprod Y_j]$. Therefore, $F(X)$ is compact.

Recall that a generating set is a set of objects, $\lbrace X_i \rbrace$, for which $[X_i,A] = 0$ for all $i$ implies that $A$ is isomorphic to the zero object. Assume that we have a set of compact objects, $\lbrace X_i \rbrace$, which is a generating set. Assume that $F(X_i)$ is compact for all $i$. Then, for each $i$,
\begin{gather*}
 [X_i,\coprod G(Y_j)] \cong \bigoplus [X_i,G(Y_j)] \cong \bigoplus [F(X_i),Y_j] \cong \\ [F(X_i),\coprod Y_j] \cong [X_i,G(\coprod Y_j)].
\end{gather*}
This morphism coincides with applying $[X_i,-]$ to the natural map $\coprod G(Y_j) \ra G(\coprod Y_j)$. Let $Z$ be the cone over $\coprod G(Y_j) \ra G(\coprod Y_j)$. $[X_i,Z] = 0$ for all $i$ by the above calculation. Thus, $Z$ is isomorphic to $0$ and $G$ commutes with coproducts. \qed

\begin{lem}
 $E$ is $\id_X$-perfect if and only if $E$ is perfect.
\end{lem}

\proof If $E$ is perfect, then, for any perfect object $F$, $E \otimes F$ is perfect. If $E$ is $\id_X$-perfect, then $\mathcal{H}om(E,-)$ commutes with coproducts. Taking global sections, we see that $[E,-]$ commutes with coproducts. \qed

\begin{lem}
 If $E$ is $f$-perfect, then the natural map,
\begin{displaymath}
 \mathcal{H}om(E,f^!G) \otimes f^*H \ra \mathcal{H}om(E,f^!(G \otimes H)),
\end{displaymath}
 is an isomorphism for any $G$ and $H$ in $D(Y)$.
\label{lem:pulloutfactor}
\end{lem}

\proof First, let us describe the origin of the aforementioned natural map. By adjunction, it is sufficient to provide a map
\begin{displaymath}
 f_*(E \otimes \mathcal{H}om(E,f^!G) \otimes f^*H) \ra G \otimes H.
\end{displaymath}
Using the projection formula, see Proposition $5.3$ of \cite{Nee96},
\begin{displaymath}
 f_*(E \otimes \mathcal{H}om(E,f^!G) \otimes f^*H) \cong f_*(E \otimes \mathcal{H}om(E,f^!G)) \otimes H.
\end{displaymath}
Using the counit of adjunction, we get a map
\begin{displaymath}
 f_*(E \otimes \mathcal{H}om(E,f^!G)) \otimes H \ra G \otimes H.
\end{displaymath}
The composition of these two maps is our desired natural map, which we call $\nu_{G,H}$. Since both 
\begin{displaymath}
 \mathcal{H}om(E,f^!-)\otimes f^*- \text{ and } \mathcal{H}om(E,f^!(-\otimes-))
\end{displaymath}
commute with coproducts, it suffices to show that $\nu_{F,F'}$ is an isomorphism for any compact $F$ and $F'$. Moreover, we only need to show that $[D,\nu_{F,F'}]$ is an isomorphism for any $D$ in $D(Y)$. We have a sequence of natural isomorphisms:
\begin{gather*}
 [D,\mathcal{H}om(E,f^!F)\otimes f^*F'] \cong [D \otimes \mathcal{H}om(f^*F',\mathcal{O}_X), \mathcal{H}om(E,f^!F)] \cong \\ [f_*(E \otimes D \otimes f^*\mathcal{H}om(F',\mathcal{O}_Y)),F] \cong [f_*(E \otimes D) \otimes \mathcal{H}om(F',\mathcal{O}_Y),F] \cong \\ [f_*(E \otimes D),F' \otimes F] \cong [D, \mathcal{H}om(E,f^!(F \otimes F'))].
\end{gather*}
It is straightforward, but tedious, to check the end result of this sequence coincides with the precomposition by $\nu_{F,F'}$. \qed

\begin{lem}
 If $E$ is $f$-perfect, then we have a natural isomorphism,
\begin{displaymath}
 f_*\mathcal{H}om(E,f^!G) \cong \mathcal{H}om(f_*E,G),
\end{displaymath}
for any $G$ in $D(Y)$.
\label{lem:sheafGrothduality}
\end{lem}

\proof From the counit, $f_*f^!G \ra G$, we get a natural map $f_*(\mathcal{H}om(E,f^!G)) \ra \mathcal{H}om(f_*E,G)$. To show that is it an isomorphism of sheaves, we need to know that it induces an isomorphism of the derived global sections over any open affine subset, $U$, of $Y$. Let $j: U \ra Y$ and $j':f^{-1}U \ra X$ be the inclusions and $f'$ the restriction of $f$ to $f^{-1}U$. It is enough to check that $(j')^*\mathcal{H}om(E,f^!-)$ is isomorphic to $\mathcal{H}om(E',(f')^!j^*-)$, where $E'$ is $(j')^*E$. Since localization is exact, we have $j^*f_* \cong f'_* (j')^*$. Take any compact $E$ from $D(X)$. We see that $j^*f_*(E \otimes -) \cong f'_*(E' \otimes (j')^*-)$. Taking adjoints, we get
\begin{displaymath}
 \mathcal{H}om(E,f^!j_*-) \cong j'_*\mathcal{H}om(E',(f')^!-).
\end{displaymath}
Consequently,
\begin{displaymath}
 (j')^*\mathcal{H}om(E,f^!-)j_*j^* \cong (j')^*j'_*\mathcal{H}om(E',(f')^!-)j^*.
\end{displaymath}
Since $(j')^*j'_*$ is the identity, we have
\begin{displaymath}
 \mathcal{H}om(E',(f')^!-)j^* \cong (j')^*\mathcal{H}om(E,f^!-)j_*j^*.
\end{displaymath}
We use the unit, $\id \ra j_*j^*$, of adjunction and reduce to checking that the resulting natural map,
\begin{displaymath}
 (j')^*\mathcal{H}om(E,f^!-) \ra (j')^*\mathcal{H}om(E,f^!-)j_*j^*,
\end{displaymath}
is an isomorphism. The cone over the unit, $\id \ra j_*j^*$, is the functor of local cohomology on $Z = Y - U$. Let us denote this by $\Gamma_Z$. It is now sufficient to check that $(j')^*\mathcal{H}om(E,f^!-)\Gamma_Z = 0$.

If $G$ is any complex acyclic off of $Z$, $\mathcal{H}om(E,f^!G) \cong \mathcal{H}om(E,f^!\mathcal{O}_Y) \otimes f^*G$ is acyclic off $f^{-1}Z$ as $f^*G$ is. Thus, after we apply $(j')^*$, we get zero. \qed

\begin{lem}
 If $E$ is $f$-perfect, then so is $\mathcal{H}om(E,f^!F)$ for any $F$ in $D_{\perf}(Y)$.
\end{lem}

\proof We need to show that $f_*(\mathcal{H}om(E,f^!F) \otimes -)$ takes perfect objects to perfect objects. By lemma \ref{lem:sheafGrothduality}, this is isomorphic to $\mathcal{H}om(f_*(E \otimes \mathcal{H}om(-,\mathcal{O}_X)),F)$. As $f_*(E\otimes \mathcal{H}om(-,\mathcal{O}_X))$ and $F$ are perfect, $\mathcal{H}om(f_*(E\otimes \mathcal{H}om(-,\mathcal{O}_X)),F)$ is also perfect.  \qed

\begin{lem}
 $E \in D(X)$ is zero if and only if $f_*(E \otimes F)$ is zero for all perfect $F$.
\label{lem:detectzerocomplexwithpushforward}
\end{lem}

\proof Clearly, if $E$ is zero, then $f_*(E \otimes G)$ is zero for any $G$ in $D(X)$. Assume that $f_*(E \otimes F)$ is zero for all perfect $F$. Then,
\begin{displaymath}
 0 = [\mathcal{O}_Y,f_*(E \otimes F)] \cong [f^*\mathcal{O}_Y, E \otimes F] \cong [\mathcal{H}om(F,\mathcal{O}_X),E]
\end{displaymath}
for any $F$. Since $\mathcal{H}om(-,\mathcal{O}_X)$ is a involution on $D_{\perf}(X)$, $[F,E] = 0$ 
for any perfect $F$. Thus, $E$ is zero. \qed

\begin{lem}
 If $E$ is $f$-perfect, the natural map
\begin{displaymath}
 \nu: E \ra \mathcal{H}om(\mathcal{H}om(E,f^!\mathcal{O}_Y),f^!\mathcal{O}_Y)
\end{displaymath}
is an isomorphism. Therefore, $\mathcal{H}om(-,f^!\mathcal{O}_Y)$ is an involution on the collection of $f$-perfect objects.
\end{lem}

\proof By lemma \ref{lem:detectzerocomplexwithpushforward}, it suffices to show that $f_*(\nu \otimes F)$ is an isomorphism for any perfect object $F$. By lemma \ref{lem:sheafGrothduality},
\begin{gather*}
 f_*\mathcal{H}om(\mathcal{H}om(E \otimes F,f^!\mathcal{O}_Y),f^!\mathcal{O}_Y) \cong \mathcal{H}om(f_*\mathcal{H}om(E \otimes F,f^!\mathcal{O}_Y),\mathcal{O}_Y) \cong \\ \mathcal{H}om(\mathcal{H}om(f_*(E \otimes F),\mathcal{O}_Y),\mathcal{O}_Y).
\end{gather*}
It is straightforward to check that the resulting map coincides with the double dualization map, over $\mathcal{O}_Y$, for $f_*(E \otimes F)$. $f_*(E \otimes F)$ is perfect so the double dualization
\begin{displaymath}
 f_*(E \otimes F) \ra \mathcal{H}om(\mathcal{H}om(f_*(E \otimes F),\mathcal{O}_Y),\mathcal{O}_Y)
\end{displaymath}
is an isomorphism. \qed

\begin{lem}
 If $E$ is $f$-perfect, then $f_*(\mathcal{H}om(E,f^!\mathcal{O}_Y) \otimes -)$ is left adjoint to $E \otimes f^*-$.
\label{lem:perftoleftadj}
\end{lem}

\proof We compute:
\begin{gather*}
 [F,E \otimes f^*G] \cong [F, \mathcal{H}om(\mathcal{H}om(E,f^!\mathcal{O}_Y),f^!\mathcal{O}_Y) \otimes f^*G] \cong \\ [F,\mathcal{H}om(\mathcal{H}om(E,f^!\mathcal{O}_Y),f^!G)] \cong [f_*(\mathcal{H}om(E,f^!\mathcal{O}_Y) \otimes F),G].
\end{gather*} \qed

We say that an object $G$ of $D(X)$ is \textbf{locally-finite} if
\begin{displaymath}
 \sum_{i \in \Z} \dim_k [F,G[i]] < \infty
\end{displaymath}
for all perfect $F$.

\begin{cor}
 Let $X$ and $Y$ be defined over a field $k$. If $E$ is $f$-perfect, $E \otimes f^*-$ takes locally-finite objects to locally-finite objects. In particular, if $X$ and $Y$ are quasi-projective over $k$, $E \otimes f^*-$ must take $D^b_{\cohc}(Y)$ to $D^b_{\cohc}(X)$ if $E$ if $f$-perfect.
\label{cor:f-perftotensorperf}
\end{cor}

\proof $f_*(\mathcal{H}om(E,f^!\mathcal{O}_Y) \otimes -)$ takes $D_{\perf}(X)$ to $D_{\perf}(Y)$. The right adjoint must therefore take locally-finite objects to locally-finite objects. If $X$ and $Y$ are quasi-projective over $k$, then locally-finite objects are exactly the objects of $D^b_{\cohc}(X)$. \qed

The following is a slight extension of proposition $1.6$ from \cite{RMS07}.

\begin{prop}
 Let $f:X \ra Y$ be a morphism of quasi-projective schemes over $k$ and $E \in D(X)$. The following are equivalent:
 \begin{itemize}
  \item $E$ is $f$-perfect.
  \item $\mathcal{H}om(E,f^!G)$ lies in $D^b_{\cohc}(X)$ for all $G$ in $D^b_{\cohc}(Y)$.
 \end{itemize}
 If $E \in D^b_{\coh}(X)$ and $f$ is proper, then the above are equivalent to:
 \begin{itemize}
  \item $E \otimes f^*G$ lies in $D^b_{\cohc}(X)$ for all $G \in D^b_{\cohc}(Y)$.
 \end{itemize}
\label{prop:charrelperf}
\end{prop}

\proof Condition one implies condition two because $\mathcal{H}om(E,f^!-)$ is right adjoint to a functor that takes perfect objects to perfect objects. Similarly, condition two implies condition one since $f_*(E \otimes -)$ is left adjoint to a functor that sends $D^b_{\cohc}(Y)$ to $D^b_{\cohc}(X)$. 

Corollary \ref{cor:f-perftotensorperf} says that condition one implies condition three. If condition three holds for $E$, it also holds for $E \otimes F$ for any perfect $F$. Using the projection formula, $f_*(E \otimes F) \otimes G$ lies in $D^b_{\cohc}(Y)$ for any $G$ in $D^b_{\cohc}(Y)$. The proof of this proposition is finished by the following lemma. \qed

\begin{lem}
 If $E \in D^b_{\coh}(Y)$ and $E \otimes G$ lies in $D^b_{\cohc}(Y)$ for all $G$ in $D^b_{\cohc}(Y)$, then $E$ is perfect.
\label{lem:perfectviatensor}
\end{lem}

\proof Assume that $E \otimes G$ lies in $D^b_{\coh,c}(Y)$ for all $G$ in $D^b_{\coh,c}(Y)$. Let $p$ be a closed point of $X$. Then $E \otimes k(p)$ must be a bounded complex of coherent sheaves. Let $E_p = E \otimes \mathcal{O}_{p,X}$. Since $E_p \otimes k(p) = E \otimes k(p)$, we see that $E_p$ has finite Tor-dimension in $\mathcal{O}_{p,X}$. Replacing $E_p$ by its minimal free resolution, we see that $E_p$ is quasi-isomorphic to bounded complex of finite rank free modules. This must be true on a neighborhood of $p$. So $E$ is locally quasi-isomorphic to a bounded complex of locally-free sheaves. \qed

We can also test for perfection by using Hom instead of $\otimes$.

\begin{lem}
 Assume $E$ lies in $D^b_{\coh}(X)$. $E$ is perfect if and only if $\mathcal{H}om(E,G)$ lies in $D^b_{\coh}(X)$ for any $G$ in $D^b_{\coh}(X)$.
\label{lem:perfviahom}
\end{lem}

\proof Clearly, if $E$ is perfect, then $\mathcal{H}om(E,G)$ is a bounded complex of coherent sheaves for any $G \in D^b_{\coh}(X)$.

Assume that $\mathcal{H}om(E,G)$ lies in $D^b_{\coh}(X)$ for any $G$ in $D^b_{\coh}(X)$. We can apply the argument from the previous lemma except we consider $\Hom(E_p,k(p))$. If $E_p$ is replaced by its minimal free resolution, then $\Hom(E_p,k(p))$ is bounded if and only if $E_p$ quasi-isomorphic to a bounded complex of free modules. $E$ must then be perfect. \qed

\begin{cor}
 If $E$ is $f$-perfect and $Y$ is projective, $E$ lies in $D^b_{\coh}(X)$.
\end{cor}

Now we can begin applying the above homological algebra to the case of integral transforms. Let $X$ and $Y$ be be quasi-projective schemes. We have the projections
\begin{center}
 \leavevmode
 \begin{xy}
  (0,5)*+{X \times Y}="a"; (-9,-5)*+{X}="b"; (9,-5)*+{Y}="c"; {\ar@{->}_{p_1} "a";"b"}; {\ar@{->}^{p_2} "a";"c"};
 \end{xy}
\end{center}
From an object, $E \in D(X \times Y)$, we can create a functor
\begin{align*}
 \Phi_E^{X \ra Y}: D(X) & \ra D(Y) \\
 \Phi_E^{X \ra Y}(F) & = p_{2*}(E \otimes p_1^*F).
\end{align*}
We shall often drop the superscript if the context is clear.

\begin{defn}
 $\Phi_E$ is called the \textbf{integral transform} associated to $E$. $E$ is called the \textbf{kernel} of $\Phi_E$.
\end{defn}

We can immediately apply the results about relative perfection to the case of an integral transform between quasi-projective schemes over a field.

\begin{lem}
 $\Phi_E$ takes $D_{\perf}(X)$ to $D_{\perf}(Y)$ if and only if $E$ is $p_2$-perfect.
\label{lem:bounded->perfect}
\end{lem}

\proof If $E$ is $p_2$-perfect, then $\Phi_E$ must take $D_{\perf}(X)$ to $D_{\perf}(Y)$. Assume that $\Phi_E$ takes $D_{\perf}(X)$ to $D_{\perf}(Y)$. To show that $E$ is $p_2$-perfect it is enough to show that $p_{2*}(E \otimes -)$ takes a generating set of perfect objects to perfect objects. $\mathcal{O}_X(1) \boxtimes \mathcal{O}_Y(1)$ is a very ample sheaf on $X \times Y$. Thus, if $p_{2*}(E\otimes -)$ takes all shifts and tensor powers of $\mathcal{O}_X(1) \boxtimes \mathcal{O}_Y(1)$ to perfect objects, then $E$ is $p_2$-perfect. Using the projection formula, we get
\begin{displaymath}
 p_{2*}(E \otimes (\mathcal{O}_X(j) \boxtimes \mathcal{O}_Y(j))[l]) \cong p_{2*}(E \otimes p_1^*\mathcal{O}_X(j)) \otimes \mathcal{O}_Y(j)[l] = \Phi_E(\mathcal{O}_X(j)) \otimes \mathcal{O}_Y(j)[l].
\end{displaymath}
$E$ is $p_2$-perfect. \qed

\begin{lem}
 Assume $X$ and $Y$ are projective and $E$ lies in $D^b_{\coh}(X \times Y)$. $\Phi_E$ takes $D^b_{\coh}(X)$ to $D^b_{\coh}(Y)$ if and only if $E$ is $p_1$-perfect.
\label{lem:boundcohcopt->perf}
\end{lem}

\proof If $E$ is $p_1$-perfect, then the proof of proposition \ref{prop:charrelperf} shows that $\Phi_E$ takes $D^b_{\coh}(X)$ to $D^b_{\coh}(Y)$.

Assume that $\Phi_E$ takes $D^b_{\coh}(X)$ to $D^b_{\coh}(Y)$. We know that $p_{2*}(E \otimes p_1^*G)$ lies in $D^b_{\coh}(Y)$ for any $G$ from $D^b_{\coh}(X)$. We want to conclude that $E \otimes p_1^*G$ lies in $D^b_{\coh}(X \times Y)$ for any $G$ from $D^b_{\coh}(X)$. First, we prove a lemma about boundness.

\begin{lem}
 Let $X$ be a quasi-projective scheme over $k$ and $G \in D(X)$. There exists an $N \geq 0$ so that $G$ is bounded, as a complex, if and only if $[\mathcal{O}_X(j),G[l]]$ is bounded in $l$ for all $0 \leq j \leq N$.
\end{lem}

\proof If $G$ is bounded, then $[\mathcal{O}_X(j),G[l]]$ is bounded in $l$ for all $j \in \Z$.

Assume that $[\mathcal{O}_X(j),G[l]]$ is bounded for all $0 \leq j \leq N$ and $l \in \Z$. It suffices to show that $[\mathcal{O}_X(-j),G[l]] = 0$ for all $j \gg 0$ and $l \not \in [k_1,k_2]$. Using the Beilinson resolution of the diagonal on $N$-dimensional projective space and pulling back to $X$, we obtain a quasi-isomorphism
\begin{displaymath}
 \mathcal{O}_X(-j) \cong \left( V^j_M \otimes_k \mathcal{O}_X \ra \cdots \ra V^j_0 \otimes_k \mathcal{O}_X(N) \right)
\end{displaymath}
for $j>0$ where $V^j_i = H^N(\P^N_k,\Omega^i(i-j-N))$ for $0 \leq i \leq N$ (see section \ref{sec:resolvediag}). Thus, if $[\mathcal{O}_X(j),G[l]]$ is bounded in $l$ for $0 \leq j \leq N$, so $[\mathcal{O}_X(-j),G[l]]$ is uniformly bounded in $l$ for all $j < 0$. \qed

To show that $E \otimes p^*_1G$ is bounded, we just need to show that $[\mathcal{O}_X(j) \boxtimes \mathcal{O}_Y(j),E \otimes p^*_1G[l]]$ is bounded in $l$ for all $0 \leq j \leq N$.
\begin{gather*}
 [\mathcal{O}_X(j) \boxtimes \mathcal{O}_Y(j),E \otimes p^*_1G[l]] \cong [p^*_2\mathcal{O}_Y(j),E \otimes p^*_1 G(-j)[l]] \\ \cong [\mathcal{O}_Y(j),p_{2*}(E \otimes p^*_1 G(-j))[l]]
\end{gather*}
which is bounded by assumption. Therefore, $E \otimes p_1^*G$ lies in $D^b_{\coh}(X \times Y)$ and we can apply proposition \ref{prop:charrelperf} to conclude that $E$ is $p_1$-perfect. \qed

Set
\begin{displaymath}
 \Psi^{X \ra Y}_E (F) = p_{2*}\mathcal{H}om(E,p_1^!F).
\end{displaymath}

\begin{lem}
 $\Psi^{Y \ra X}_E$ is right adjoint to $\Phi^{X \ra Y}_E$.
\end{lem}

\proof This is a simple application of standard adjunctions:
\begin{gather*}
 [\Phi^{X \ra Y}_E(F),G] = [p_{2*}(E \otimes p_1^*F),G] \cong [E \otimes p_1^*F,p_2^!G] \cong \\ [p_1^*F,\mathcal{H}om(E,p_2^!G)] = [F,p_{1*}\mathcal{H}om(E,p_2^!F)].
\end{gather*} \qed

By uniqueness, $\Psi_E^{Y \ra X}|_{D^b_{\cohc}(X)}$ is the right pseudo-adjoint to $\Phi_E^{X \ra Y}|_{D_{\perf}(X)}$ if $E$ is $p_2$-perfect. From lemma \ref{lem:pulloutfactor} and proposition \ref{prop:charrelperf}, we get the following statement.
\begin{lem}
 If $E$ is $p_2$-perfect, then $\Psi_E^{Y \ra X}$ is isomorphic to $\Phi^{Y \ra X}_{\mathcal{H}om(E,p_2^!\mathcal{O}_Y)}$ and takes $D^b_{\cohc}(Y)$ to $D^b_{\cohc}(X)$. 
\end{lem}

The following tells us when an integral transform cannot have a left adjoint.

\begin{lem}
 If $\Phi_E$ does not take $D^b_{\cohc}(X)$ to $D^b_{\cohc}(Y)$, then it does not possess a left adjoint.
\end{lem}

\proof Assume such a left adjoint exists. It would have to take perfect objects to perfect objects since $\Phi_E$ commutes with coproducts. By duality, $\Phi_E$ would then have to take $D^b_{\cohc}(X)$ to $D^b_{\cohc}(Y)$. \qed

\begin{lem}
 Assume that $X$ and $Y$ are projective and $E$ lies in $D^b_{\coh}(X \times Y)$. If $\Phi^{X \ra Y}_E$ takes $D^b_{\cohc}(X)$ to $D^b_{\cohc}(Y)$, then $\Phi^{Y \ra X}_{\mathcal{H}om(E,p_1^!\mathcal{O}_X)}$ is left adjoint to $\Phi^{X \ra Y}_E$.
\end{lem}

\proof From lemma \ref{lem:boundcohcopt->perf}, we know that $E$ is $p_1$-perfect. Take $F$ from $D(Y)$ and $G$ from $D(X)$. We apply lemma \ref{lem:perftoleftadj} in the following sequence of natural isomorphisms:
\begin{gather*}
 [F,\Phi^{X \ra Y}_E(G)] \cong [p_2^*F,E \otimes p_1^*G] \cong \\ [p_{1*}(\mathcal{H}om(E,p_1^!\mathcal{O}_X) \otimes p^*_2 F), G] = [\Phi^{Y \ra X}_{\mathcal{H}om(E,p_1^!\mathcal{O}_X)}(F),G].
\end{gather*}
\qed

\section{Totalizing complexes in triangulated categories}
\label{sec:tot}
 
The majority of the ideas and results of this section are due to Orlov \cite{Orl97}, see also \cite{Kaw04}. Let $\mathcal{T}$ be a triangulated category.

\begin{defn}
 A \textbf{complex over $\mathcal{T}$} is a diagram
\begin{center}
 \leavevmode
 \begin{xy}
  (-50,7.5)*+{A_s}="a"; (-30,7.5)*+{A_{s+1}}="b"; (-10,7.5)*+{A_{s+2}}="c"; (10,7.5)*+{\cdots}="d"; (30,7.5)*+{A_{-1}}="e"; (50,7.5)*+{A_0}="f"; {\ar@{->}^{d_s} "a";"b"};  {\ar@{->}^{d_{s+1}} "b";"c"}; {\ar@{->}^{d_{s+2}} "c";"d"}; {\ar@{->}^{d_{-2}} "d";"e"}; {\ar@{->}^{d_{-1}} "e";"f"};
 \end{xy}
\end{center}
 with $A_i \in \mathcal{T}$, $d_i \in [A_i,A_{i+1}]$, and $d_{i+1} \circ d_i = 0$.
 A morphism $f$ between two complexes, $(A,d)$ and $(A',d')$, over $\mathcal{T}$ is a collection of morphisms, $f_i: A_i \ra A_i'$, rendering the diagram
\begin{center}
 \leavevmode
 \begin{xy}
  (-50,7.5)*+{A_s}="a"; (-30,7.5)*+{A_{s+1}}="b"; (-10,7.5)*+{A_{s+2}}="c"; (10,7.5)*+{\cdots}="d"; (30,7.5)*+{A_{-1}}="e"; (50,7.5)*+{A_0}="f";
  (-50,-7.5)*+{A_s'}="A"; (-30,-7.5)*{A_{s+1}'}="B"; (-10,-7.5)*+{A_{s+2}'}="C"; (10,-7.5)*+{\cdots}="D"; (30,-7.5)*+{A_{-1}'}="E"; (50,-7.5)*+{A_0'}="F";
  {\ar@{->}^{f_s} "a";"A"}; {\ar@{->}^{d_s} "a";"b"}; {\ar@{->}^{d_s'} "A";"B"}; {\ar@{->}^{f_{s+1}} "b";"B"}; {\ar@{->}^{d_{s+1}} "b";"c"}; {\ar@{->}^{d_{s+1}'} "B";"C"}; {\ar@{->}^{d_{s+2}'} "C";"D"}; {\ar@{->}^{d_{-2}'} "D";"E"}; {\ar@{->}^{d_{-1}'} "E";"F"}; {\ar@{->}^{d_{s+2}} "c";"d"}; {\ar@{->}^{d_{-2}} "d";"e"}; {\ar@{->}^{d_{-1}} "e";"f"}; {\ar@{->}^{f_{s+2}} "c";"C"}; {\ar@{->}^{f_{-1}} "e";"E"}; {\ar@{->}^{f_0} "f";"F"};
 \end{xy}
\end{center}
commutative.
\end{defn}

\begin{defn}
 A \textbf{right convolution} of a complex $(A,d)$ over $\mathcal{T}$ is an object, $X \in \mathcal{T}$, and a map, $r: A_0 \ra X$, such that there exists a diagram
\begin{center}
 \leavevmode
 \begin{xy}
  (-55,7.5)*+{A_s}="a"; (-35,7.5)*+{A_{s+1}}="b"; (-15,7.5)*+{A_{s+2}}="c"; (5,7.5)*+{\cdots}="d"; (25,7.5)*+{A_{-1}}="e"; (45,7.5)*+{A_0}="f";
  (-45,-7.5)*+{A_s=B_s}="A"; (-25,-7.5)*{B_{s+1}}="B"; (-5,-7.5)*+{B_{s+2}}="C"; (15,-7.5)*+{\cdots}="D"; (35,-7.5)*+{B_{-1}}="E"; (55,-7.5)*+{B_0=X}="F";
  (-45,0)*+{\circlearrowleft}; (-35,0)*+{\vartriangle}; (-25,0)*+{\circlearrowleft}; (-15,0)*+{\vartriangle}; (-5,0)*+{\circlearrowleft}; (25,0)*+{\vartriangle}; (35,0)*+{\circlearrowleft}; (45,0)*+{\vartriangle};
  {\ar@{->}_= "a";"A"}; {\ar@{->}^{d_s} "a";"b"}; {\ar@{-->} "B";"A"}; {\ar@{->} "A";"b"}; {\ar@{->} "b";"B"}; {\ar@{->}^{d_{s+1}} "b";"c"}; {\ar@{-->} "C";"B"}; {\ar@{-->} "D";"C"}; {\ar@{-->} "E";"D"}; {\ar@{-->} "F";"E"}; {\ar@{->}^{d_{s+2}} "c";"d"}; {\ar@{->}^{d_{-2}} "d";"e"}; {\ar@{->}^{d_{-1}} "e";"f"}; {\ar@{->} "B";"c"}; {\ar@{->} "C";"d"}; {\ar@{->} "D";"e"}; {\ar@{->} "E";"f"}; {\ar@{->} "c";"C"}; {\ar@{->} "e";"E"}; {\ar@{->}^{r} "f";"F"};
 \end{xy}
\end{center}
where faces denoted by $\circlearrowleft$ are commutative and faces denoted by $\vartriangle$ are exact triangles in $\mathcal{T}$.
 A \textbf{left convolution} of a complex $(A,d)$ over $\mathcal{T}$ is an object, $Y \in \mathcal{T}$, and a map, $l: Y \ra A_s$, such that there exists a diagram
\begin{center}
 \leavevmode
 \begin{xy}
  (-45,7.5)*+{A_s}="a"; (-25,7.5)*+{A_{s+1}}="b"; (-5,7.5)*+{A_{s+2}}="c"; (15,7.5)*+{\cdots}="d"; (35,7.5)*+{A_{-1}}="e"; (55,7.5)*+{A_0}="f";
  (-55,-7.5)*+{Y=B_s}="A"; (-35,-7.5)*{B_{s+1}}="B"; (-15,-7.5)*+{B_{s+2}}="C"; (5,-7.5)*+{\cdots}="D"; (25,-7.5)*+{B_1}="E"; (45,-7.5)*+{B_0=A_0}="F";
  (-45,0)*+{\vartriangle}; (-35,0)*+{\circlearrowleft}; (-25,0)*+{\vartriangle}; (-15,0)*+{\circlearrowleft}; (-5,0)*+{\vartriangle}; (25,0)*+{\circlearrowleft}; (35,0)*+{\vartriangle}; (45,0)*+{\circlearrowleft};
  {\ar@{->}^l "A";"a"}; {\ar@{->}^{d_s} "a";"b"}; {\ar@{-->} "B";"A"}; {\ar@{->} "a";"B"}; {\ar@{->} "B";"b"}; {\ar@{->}^{d_{s+1}} "b";"c"}; {\ar@{-->} "C";"B"}; {\ar@{-->} "D";"C"}; {\ar@{-->} "E";"D"}; {\ar@{-->} "F";"E"}; {\ar@{->}^{d_{s+2}} "c";"d"}; {\ar@{->}^{d_{-2}} "d";"e"}; {\ar@{->}^{d_{-1}} "e";"f"}; {\ar@{->} "b";"C"}; {\ar@{->} "c";"D"}; {\ar@{->} "d";"E"}; {\ar@{->} "e";"F"}; {\ar@{->} "C";"c"}; {\ar@{->} "E";"e"}; {\ar@{->}_= "F";"f"};
 \end{xy}
\end{center}
where faces denoted by $\circlearrowleft$ are commutative and faces denoted by $\vartriangle$ are exact triangles in $\mathcal{T}$.
\end{defn}

\begin{lem}
 Let $(A,d)$ be a complex over $\mathcal{T}$. Assume that $[A_l,A_m[j]]$ is zero for all $l < m$ and $j < 0$. Then, there exists a right convolution $X$ of $(A,d)$.
\label{lem:rightconvolve}
\end{lem}

\proof We proceed by induction on the length of the complex. Assume that the lemma is true for any complex that satisfies the hypotheses of the lemma and has length less than $-s$. Form the triangle
\begin{center}
 \leavevmode
 \begin{xy}
  (-10,7)*+{A_s}="a"; (10,7)*+{A_{s+1}}="b"; (0,-7)*+{Y}="c"; {\ar@{-->} "c";"a"}; {\ar@{->}^{d_s} "a";"b"}; {\ar@{->} "b";"c"};
 \end{xy}
\end{center}
and consider the long exact sequence
\begin{center}
 \leavevmode
 \begin{xy}
  (-50,0)*+{\cdots}="a"; (-30,0)*+{[Y,A_{s+2}]}="b"; (0,0)*+{[A_{s+1},A_{s+2}]}="c"; (30,0)*+{[A_s,A_{s+2}]}="d"; (50,0)*+{\cdots}="e"; {\ar@{->} "a"; "b"}; {\ar@{->} "b"; "c"}; {\ar@{->} "c"; "d"}; {\ar@{->} "d"; "e"};
 \end{xy}
\end{center}
resulting from applying $[-,A_{s+2}]$. Since $d_{s+1} \circ d_s = 0$, there is a map $e: Y \ra A_{s+2}$ making
\begin{center}
 \leavevmode
 \begin{xy}
  (-10,7)*+{A_{s+1}}="b"; (-10,-7)*+{Y}="c"; (10,7)*+{A_{s+2}}="d"; {\ar@{->} "b";"c"}; {\ar@{->}^{d_{s+1}} "b";"d"}; {\ar@{->}^{e} "c";"d"};
 \end{xy}
\end{center}
commute. Now consider the long exact sequence
\begin{center}
 \leavevmode
 \begin{xy}
  (-55,0)*+{\cdots}="a"; (-30,0)*+{[A_s[1],A_{s+3}]}="b"; (0,0)*+{[Y,A_{s+3}]}="c"; (30,0)*+{[A_{s+1},A_{s+3}]}="d"; (55,0)*+{\cdots}="e"; {\ar@{->} "a"; "b"}; {\ar@{->} "b"; "c"}; {\ar@{->} "c"; "d"}; {\ar@{->} "d"; "e"};
 \end{xy}
\end{center}
resulting from applying $[-,A_{s+3}]$. $d_{s+2} \circ e \in [Y,A_{s+3}]$ maps to $d_{s+2} \circ d_{s+1} = 0$ and $[A_s[1],A_{s+3}] = 0$ by assumption. Thus, $d_{s+2} \circ e = 0$. By examining the long exact sequence resulting from applying $[-,A_l[j]]$, we see that $[Y,A_l[j]]$ is zero for $s+2 \leq l \leq 0$ and $j < 0$. Thus,
\begin{center}
 \leavevmode
 \begin{xy}
  (-50,7.5)*+{Y}="a"; (-30,7.5)*+{A_{s+2}}="b"; (-10,7.5)*+{A_{s+3}}="c"; (10,7.5)*+{\cdots}="d"; (30,7.5)*+{A_{-1}}="e"; (50,7.5)*+{A_0}="f"; {\ar@{->}^{e} "a";"b"};  {\ar@{->}^{d_{s+2}} "b";"c"}; {\ar@{->}^{d_{s+3}} "c";"d"}; {\ar@{->}^{d_{-2}} "d";"e"}; {\ar@{->}^{d_{-1}} "e";"f"};
 \end{xy}
\end{center}
is a complex over $\mathcal{T}$ satisfying the hypotheses of the lemma and having length less than $-s$. It possesses a right convolution, which is also a right convolution of the original complex. \qed

\begin{lem}
 Let $(A,d)$ and $(A',d')$ be complexes over $\mathcal{T}$ so that $[A_l,A_m[j]]$ is zero for $l < m$ and $j<0$ and $[A_l',A_m'[j]]$ is zero for $l < m$ and $j < 0$. Assume that $[A_l,A_m'[j]]$ is zero for $l < m$ and $j<0$. Then, for any right convolutions $(X,r)$ of $(A,d)$ and $(X',r')$ of $(A',d')$ and any map $h: X' \ra Z'$ there exists a (non-canonical) morphism $g: X \ra Z'$ making the diagram
\begin{center}
 \leavevmode
 \begin{xy}
  (-20,10)*+{A_0}="a"; (0,10)*+{X}="b"; (20,10)*+{X}="c"; (-20,-10)*+{A_0'}="A"; (0,-10)*+{X'}="B"; (20,-10)*+{Z'}="C"; {\ar@{->}_{f_0} "a";"A"}; {\ar@{->}^r "a";"b"}; {\ar@{->}^= "b";"c"}; {\ar@{->}^{r'} "A";"B"}; {\ar@{->}^h "B";"C"}; {\ar@{->}^g "c";"C"};
 \end{xy}
\end{center}
commute. If, in addition, $[A_l,Z'[j]]$ is zero for all $l,j <0$, then this morphism is unique.
\label{lem:rightconvolvemorphisms}
\end{lem}

\proof We again proceed by induction. Assume the result is true for morphisms satisfying the hypotheses of the lemma and having length less than $-s$. Form the triangles
\begin{center}
 \leavevmode
 \begin{xy}
  (-30,7)*+{A_s}="a"; (-10,7)*+{A_{s+1}}="b"; (-20,-7)*+{Y}="c"; {\ar@{-->} "c";"a"}; {\ar@{->}^{d_s} "a";"b"}; {\ar@{->} "b";"c"}; (10,7)*+{A_s'}="a"; (30,7)*+{A_{s+1}'}="b"; (20,-7)*+{Y'}="c"; {\ar@{-->} "c";"a"}; {\ar@{->}^{d_s'} "a";"b"}; {\ar@{->} "b";"c"};
 \end{xy}
\end{center}
and let $e:Y \ra A_{s+2}$ and $e':Y' \ra A_{s+2}'$ be the maps constructed as in the previous lemma. There is a $y: Y \ra Y'$ making
\begin{center}
 \leavevmode
 \begin{xy}
  (-30,10)*+{A_s}="a"; (-10,10)*+{A_{s+1}}="b"; (10,10)*+{Y}="c"; (30,10)*+{A_s[1]}="d"; (-30,-10)*+{A_s'}="A"; (-10,-10)*+{A_{s+1}'}="B"; (10,-10)*+{Y'}="C"; (30,-10)*+{A_s'[1]}="D"; {\ar@{->}^{d_s} "a";"b"}; {\ar@{->} "b";"c"}; {\ar@{->} "c";"d"}; {\ar@{->}^{d_s'} "A";"B"}; {\ar@{->} "B";"C"}; {\ar@{->} "C";"D"}; {\ar@{->}_{f_s} "a";"A"}; {\ar@{->}_{f_{s+1}} "b";"B"}; {\ar@{->}_{f_s[1]} "d";"D"}; {\ar@{->}_y "c";"C"};
 \end{xy}
\end{center}
commutative. We have a diagram
\begin{center}
 \leavevmode
 \begin{xy}
  (-20,15)*+{A_{s+1}}="a"; (0,7)*+{Y}="b"; (20,15)*+{A_{s+2}}="c"; (-20,-15)*+{A_{s+1}'}="A"; (0,-7)*+{Y'}="B"; (20,-15)*+{A_{s+2}'}="C"; {\ar@{->} "a";"b"}; {\ar@{->}_{e} "b";"c"}; {\ar@{->}^{d_{s+1}} "a";"c"}; {\ar@{->} "A";"B"}; {\ar@{->}^{e'} "B";"C"}; {\ar@{->}_{d_{s+1}'} "A";"C"}; {\ar@{->}_{f_{s+1}} "a";"A"}; {\ar@{->}_{y} "b";"B"}; {\ar@{->}_{f_{s+2}} "c";"C"};
 \end{xy}
\end{center}
with all but the right parallelogram commutative. Consider the long exact sequence
\begin{center}
 \leavevmode
 \begin{xy}
  (-50,0)*+{\cdots}="a"; (-30,0)*+{[A_s[1],A_{s+2}']}="b"; (0,0)*+{[Y,A_{s+2}']}="c"; (30,0)*+{[A_{s+1},A_{s+2}']}="d"; (50,0)*+{\cdots}="e"; {\ar@{->} "a"; "b"}; {\ar@{->} "b"; "c"}; {\ar@{->} "c"; "d"}; {\ar@{->} "d"; "e"};
 \end{xy}
\end{center}
and note that $e' \circ y$ and $f_{s+2} \circ e$ both map to $d_{s+1}' \circ f_{s+1} = f_{s+2} \circ d_{s+1}$ and that $[A_s[1],A_{s+2}']$ is zero by assumption. Thus, $e' \circ y = f_{s+2} \circ e$ and we reduce to a map of complexes satisfying the hypotheses of the lemma and of length less than $-s$.

Let $(X,r)$ and $(X',r')$ be convolutions of $(A,d)$ and $(A',d')$, respectively. And, let $h: X' \ra Z'$ be a morphism. We use the notation from the definition of a convolution. Let $C$ be the cone over the map, $h \circ r': X' \ra Z'$. We have a morphism, $B_{-1} \ra C[-1]$, which makes the diagram
\begin{center}
 \leavevmode
 \begin{xy}
  (-30,10)*+{B_{-1}}="a"; (-10,10)*+{A_0}="b"; (10,10)*+{X}="c"; (30,10)*+{B_{-1}[1]}="d"; (-30,-10)*+{C[-1]}="A"; (-10,-10)*+{A_0'}="B"; (10,-10)*+{Z'}="C"; (30,-10)*+{C}="D"; {\ar@{->} "a";"b"}; {\ar@{->}^r "b";"c"}; {\ar@{->} "c";"d"}; {\ar@{->} "A";"B"}; {\ar@{->}^{h \circ r'} "B";"C"}; {\ar@{->} "C";"D"}; {\ar@{->} "a";"A"}; {\ar@{->}_{f_0} "b";"B"}; {\ar@{->} "d";"D"}; {\ar@{->}_g "c";"C"};
 \end{xy}
\end{center}
commutative. If $[A_l,Z'[j]]$ is zero for $l \leq 0$ and $j < 0$, then $[B_{-1},Z'[j]]$ is zero for $j < 0$. The long exact sequence
\begin{center}
 \leavevmode
 \begin{xy}
  (-50,0)*+{\cdots}="a"; (-30,0)*+{[B_{-1}[1],Z']}="b"; (0,0)*+{[X,Z']}="c"; (30,0)*+{[A_0,Z']}="d"; (50,0)*+{\cdots}="e"; {\ar@{->} "a"; "b"}; {\ar@{->} "b"; "c"}; {\ar@{->} "c"; "d"}; {\ar@{->} "d"; "e"};
 \end{xy}
\end{center}
shows that, if $g \circ r = h \circ r' \circ f_0$ and $g' \circ r = h \circ r' \circ f_0$, then $g=g'$. Thus, the map $g: X \ra Z'$ is unique. \qed

\begin{cor}
  Let $(A,d)$ be a complex over $\mathcal{T}$. Assume that $[A_l,A_m[j]]$ is zero for all $l < m$ and $j < 0$. All right convolutions are (non-canonically) isomorphic. If we assume, in addition, that $[A_t,X[j]]$ is zero for all $l$, $j \leq 0$, and some right convolution $X$, then all right convolutions of $(A,d)$ are canonically isomorphic.
\label{cor:rightconvolve}
\end{cor}

\proof Apply lemma \ref{lem:rightconvolvemorphisms} to the identity map between two complexes and note that, if all $f_i$ are isomorphisms, the resulting morphism between the convolution is an isomorphism. \qed

We also have duals of these results.

\begin{lem}
 Let $(A,d)$ be a complex over $\mathcal{T}$ with $[A_l,A_m[j]]$ is zero for $l<m$ and $j<0$. There exists a left convolution $(Y,l)$ of $(A,d)$ and all left convolutions are (non-canonically) isomorphic. If, in addition, $[Y,A_m[j]]$ is zero for $j<0$, $Y$ is unique up to a unique isomorphism.
 \label{lem:leftconvolve}
\end{lem}

\begin{lem}
 Let $(A,d)$ and $(A',d')$ be complexes over $\mathcal{T}$ so that $[A_l,A_m[j]]$ is zero for $l < m$ and $j<0$ and $[A_l',A_m'[j]]$ is zero for $l < m$ and $j < 0$. Assume that $[A_l,A_m'[j]]=0$ for $l < m$ and $j<0$. Then, for any left convolutions $(Y,l)$ of $(A,d)$ and $(Y',l')$ of $(A',d')$ and any morphism $h: Z \ra Y$, there exists a (non-canonical) morphism $g: Z \ra Y'$ making the diagram
\begin{center}
 \leavevmode
 \begin{xy}
  (-20,10)*+{Z}="a"; (0,10)*+{Y}="b"; (20,10)*+{A_s}="c"; (-20,-10)*+{Y'}="A"; (0,-10)*+{Y'}="B"; (20,-10)*+{A_s'}="C"; {\ar@{->}_{g} "a";"A"}; {\ar@{->}^h "a";"b"}; {\ar@{->}^l "b";"c"}; {\ar@{->}^= "A";"B"}; {\ar@{->}^{l'} "B";"C"}; {\ar@{->}^{f_s} "c";"C"};
 \end{xy}
\end{center}
commute. If, in addition, $[Z,A_l'[j]]$ is zero for all $l,j <0$, then this morphism is unique.
\label{lem:leftconvolvemorphisms}
\end{lem}

\begin{eg}
 If $\mathcal{T}=D(X)$ and each $A_i$ is a quasi-coherent sheaf placed in degree zero, then we can convolve and the convolution is simply the complex itself as an object of $D(X)$.
\end{eg}

We can also use these results to totalize unbounded complexes over $\mathcal{T}$. This is done by taking a homotopy colimit of the convolutions of the brutal truncations. Assume we are given a bounded above complex
\begin{center}
 \leavevmode
 \begin{xy}
  (-50,7.5)*+{\cdots}="a"; (-30,7.5)*+{A_{s+1}}="b"; (-10,7.5)*+{A_{s+2}}="c"; (10,7.5)*+{\cdots}="d"; (30,7.5)*+{A_{-1}}="e"; (50,7.5)*+{A_0}="f"; {\ar@{->}^{d_s} "a";"b"};  {\ar@{->}^{d_{s+1}} "b";"c"}; {\ar@{->}^{d_{s+2}} "c";"d"}; {\ar@{->}^{d_{-2}} "d";"e"}; {\ar@{->}^{d_{-1}} "e";"f"};
 \end{xy}
\end{center}
over a triangulated category $\mathcal{T}$ possessing small coproducts. The brutal $s$th truncation of $(A,d)$ is the complex 
\begin{center}
 \leavevmode
 \begin{xy}
  (-50,7.5)*+{A_s}="a"; (-30,7.5)*+{A_{s+1}}="b"; (-10,7.5)*+{A_{s+2}}="c"; (10,7.5)*+{\cdots}="d"; (30,7.5)*+{A_{-1}}="e"; (50,7.5)*+{A_0}="f"; {\ar@{->}^{d_s} "a";"b"};  {\ar@{->}^{d_{s+1}} "b";"c"}; {\ar@{->}^{d_{s+2}} "c";"d"}; {\ar@{->}^{d_{-2}} "d";"e"}; {\ar@{->}^{d_{-1}} "e";"f"};
 \end{xy}
\end{center}
It is denoted by $(\sigma_{\geq s}A,d)$. Assume that $[A_l,A_m[j]]$ is zero for $l < m$ and $j<0$. Then, we can convolve $(\sigma_{\geq s}A,d)$. Denote the convolution by $X_s$. Using the obvious morphism $\sigma_{\geq s} A \ra \sigma_{\geq s-1} A$, we get a morphism $X_s \ra X_{s-1}$. Then, we set $\text{Tot}(A,d) = \hocolim X_s$. $\text{Tot}(A,d)$ is determined up to (a non-canonical) isomorphism by $(A,d)$. 

\begin{eg}
 Let
\begin{center}
 \leavevmode
 \begin{xy}
  (-50,7.5)*+{\cdots}="a"; (-30,7.5)*+{A_{s+1}}="b"; (-10,7.5)*+{A_{s+2}}="c"; (10,7.5)*+{\cdots}="d"; (30,7.5)*+{A_{-1}}="e"; (50,7.5)*+{A_0}="f"; {\ar@{->}^{d_s} "a";"b"};  {\ar@{->}^{d_{s+1}} "b";"c"}; {\ar@{->}^{d_{s+2}} "c";"d"}; {\ar@{->}^{d_{-2}} "d";"e"}; {\ar@{->}^{d_{-1}} "e";"f"};
 \end{xy}
\end{center}
 be a complex over $D(X)$ with each $A_i$ a quasi-coherent sheaf in degree zero. Then, from the construction of the convolution, the convolution of $\sigma_{\geq s}A$ is just the complex itself as an object in $D(X)$. We shall make no distinction between the two in notation. Let $A$ denote the complex as an object of $D(X)$. Note that $A$ is the colimit of the $\sigma_{\geq s}A$ in the category of chain complexes. Therefore, there is a short exact sequence
\begin{center}
 \leavevmode
 \begin{xy}
  (-40,0)*+{0}="a"; (-20,0)*+{\bigoplus \sigma_{\geq s}A}="b"; (0,0)*+{\bigoplus \sigma_{\geq s}A}="c"; (20,0)*+{A}="d"; (35,0)*+{0}="e"; {\ar@{->} "a";"b"};  {\ar@{->} "b";"c"}; {\ar@{->} "c";"d"}; {\ar@{->} "d";"e"};
 \end{xy}
\end{center}
 of chain complexes. This induces an exact triangle in $D(X)$ and the map $\bigoplus \sigma_{\geq s}A \ra \bigoplus \sigma_{\geq s}A$ is the same as in the definition of the homotopy colimit. Thus, $\text{Tot}(A,d) \cong A$.
\end{eg}

If $\text{Tot}(A,d)$ has bounded and coherent cohomology, then there are no phantom maps $\text{Tot}(A,d) \ra \text{Tot}(A,d)$. Thus, the homotopy colimit is unique up to a unique isomorphism, given the uniqueness of the convolutions of the $\sigma_{\geq s}A$. 
 
\section{Ample sequences}
\label{sec:ample}
 
Let $\mathcal{A}$ be a $k$-linear abelian category. The following definition is due to Orlov.

\begin{defn}
 Let $\lbrace L_i \rbrace_{i \in \Z}$ be a sequence of objects in $\mathcal{A}$. We say that $\lbrace L_i \rbrace$ is an \textbf{ample sequence} if, for any object $A \in \mathcal{A}$, there exists an $N \in \Z$ so that for $i < N$ the following conditions hold:
\begin{itemize}
 \item The canonical map $[L_i,A]\otimes_k L_i \ra A$ is an epimorphism.
 \item $[L_i,A[j]]$ is zero for any $j \not = 0$.
 \item $[A,L_i]$ is zero.
\end{itemize}
\end{defn}

\begin{lem}
 Let $X$ be a quasi-projective scheme over $k$. Assume that, on $X$, we have an ample line bundle, $L$, with $L^{\otimes l} \not = \mathcal{O}_X$ for all $l \not= 0$ and $H^k(X,L^{\otimes l})$ is zero for $k>0$ and $l \gg 0$. Then, $\{L^{\otimes i}\}_{i \in \Z}$ form an ample sequence for $\Coh(X)$.
\end{lem}

\proof The first condition is a classical result of Serre. The second condition can be restated as the vanishing of $[L^{\otimes i},A[j]] = H^j(X,A \otimes L^{\otimes -i})$ for $i \ll 0$. This follows from a standard argument as in \cite{Har77}.

If $H^0(X,L^{\otimes i-j})$ is nonzero for $i-j < 0$, then $L^{\otimes j-i}$ is trivial as it is effective and its inverse has a section. Thus, $[L^{\otimes j},L^{\otimes i}]$ is zero for $j < i$. We can find a surjective map $[L^{\otimes j},A]\otimes L^{\otimes j} \ra A$ and apply $[-,L^{\otimes i}]$ to get an injection $[A,L^{\otimes i}] \ra [L^{\otimes j},A] \otimes [L^{\otimes j},L^{\otimes i}]$ for $j < i$. Consequently, $[A,L^{\otimes i}]$ is zero for $i \ll 0$. \qed

\begin{rmk}
 From the proof of the previous lemma, we see that we can replace condition three with the condition that $[L_j,L_i]$ is zero for fixed $j$ and $i \ll 0$.
\end{rmk}

\begin{lem}
 Let $\{L_i\}$ be an ample sequence in an abelian category $\mathcal{A}$. If $X$ is an object in $D^b(\mathcal{A})$ such that $\Hom(L_i,X[j]) = 0$ for all $i \ll 0$ and all $j$, then $X$ is isomorphic to zero.
\label{lem:leftspan}
\end{lem}

\proof Let $m_0$ be the minimal degree $l$ for which $H^l(X)$ is nonzero. Then, for $i \ll 0$, there is a surjection $[L_i,H^{m_0}(X)] \otimes L_i \ra H^{m_0}(X)$ and $[L_i[-m_0],X] \cong [L_i,H^{m_0}(X)]$. Thus, either $X$ is zero, or $[L_i,H^{m_0}] \not= 0$. \qed

\begin{defn}
 Let $\{P_i\}$ be an ample sequence of perfect coherent sheaves in $\Coh(X)$. We shall commonly call such a collection a \textbf{perfect ample sequence}.
\end{defn}

\begin{lem}
 Let $\{P_i\}$ be a perfect ample sequence for $X$. Then the $\{P_i\}$ generates $D_{\perf}(X)$ up to idempotent splittings.
\label{lem:perftoperf}
\end{lem}

\proof Since the set of all locally-free coherent sheaves generates $D_{\perf}(X)$ as a triangulated category, it is sufficient to show we can get any locally-free coherent sheaf from the $\{P_i\}$ by finite iteration of the operations of forming cones, forming finite direct sums, and forming direct summands. Let $V$ be any locally-free coherent sheaf. Using the first property of an ample sequence, we see that can find a resolution
\begin{displaymath}
 P_{n_i}^{\oplus m_i} \ra P_{n_{i-1}}^{\oplus m_{i-1}} \ra \cdots \ra P_{n_1}^{\oplus m_1} \ra V \ra 0.
\end{displaymath}
If $X$ is of dimension $d$, then the only map between $\ker(P_{n_i}^{\oplus m_i} \ra P_{n_{i-1}}^{\oplus m_{i-1}})[n+1]$ and $V$ if $n \geq d$ is the zero map. Thus, $V$ is a direct summand of
\begin{displaymath}
 P_{n_i}^{\oplus m_i} \ra P_{n_{i-1}}^{\oplus m_{i-1}} \ra \cdots \ra P_{n_1}^{\oplus m_1}.
\end{displaymath}
\qed

\begin{lem}
 If $Q$ is a perfect object such that $\Hom(Q,P_i[j])$ is zero for all $i \ll 0$ and all $j$. Then, $Q$ is isomorphic to zero.
\label{lem:rightspan}
\end{lem}

\proof We use the Rouquier functor, $R$, for inclusion of $D_{\perf}(X)$ into $D(X)$, see \cite{Bal09a}. Then,
\begin{displaymath}
 \Hom(Q,P_i[j])^* \cong \Hom(P_i[j],RQ)
\end{displaymath}
and $RQ$ is isomorphic to zero. Applying the duality again, we have
\begin{displaymath}
 \Hom(Q,Q)^* \cong \Hom(Q,RQ) = 0.
\end{displaymath}
Consequently, $Q$ is also isomorphic to zero. \qed

\begin{defn}
 Given a triangulated category $\mathcal{T}$, we say that a collection of objects $\lbrace S_i \rbrace_{i \in I}$ is a \textbf{spanning class} if satisfies the two following conditions:
\begin{enumerate}
 \item $[S_i,A[j]]=0$ for all $i,j$ implies that $A \cong 0$.
 \item $[A,S_i[j]]=0$ for all $i,j$ implies that $A \cong 0$.
\end{enumerate}
\end{defn}

\begin{cor}
 $\{P_i\}$ forms a spanning class for $D_{\perf}(X)$.
\end{cor}

\begin{lem}
 Let $\mathcal{T}$ be a triangulated category possesing a spanning class $\{S_i\}$. Let $F: \mathcal{T} \ra \mathcal{S}$ be an exact functor to another triangulated category possessing a left and right adjoint. If the maps
\begin{displaymath}
 \Hom(S_i,S_j[k]) \overset{\sim}{\longrightarrow} \Hom(F(S_i),F(S_j)[k])
\end{displaymath}
are isomorphisms for all $i,j,k$. Then $F$ is full and faithful.
\end{lem}

\proof Let $\leftexp{\vee}{F}$ denote the left adjoint and $F^{\vee}$ the right adjoint. Take the unit of adjunction applied to $S_i$, $f_i: S_i \ra F^{\vee} F S_i$, and form a triangle
\begin{center}
 \leavevmode
 \begin{xy}
  (-10,10)*+{S_i}="a"; (10,10)*+{F^{\vee} F S_i}="b"; (0,-5)*+{C_i}="c"; {\ar@{->}^{f_i} "a";"b"}; {\ar@{->} "b";"c"}; {\ar@{->}^{[1]} "c";"a"}
 \end{xy}
\end{center}
For all $j$ and $k$, $\Hom(S_j,C_i[k]) = 0$. From the definition of a spanning class, we conclude that $C_i \cong 0$.

Consider the counit of adjunction, $g_Q: \leftexp{\vee}{F} F Q \ra Q$, and form a triangle
\begin{center}
 \leavevmode
 \begin{xy}
  (-10,10)*+{\leftexp{\vee}{F} F Q}="a"; (10,10)*+{X}="b"; (0,-5)*+{C_Q}="c"; {\ar@{->}^{g_Q} "a";"b"}; {\ar@{->} "b";"c"}; {\ar@{->}^{[1]} "c";"a"}
 \end{xy}
\end{center}
The isomorphisms
\begin{displaymath}
 \Hom(Q,S_i[k]) \overset{\sim}{\longrightarrow} \Hom(Q,F^{\vee} F S_i[k]) \cong \Hom(\leftexp{\vee}{F} F Q,S_i[k])
\end{displaymath}
imply that $\Hom(C_Q,S_i[k]) = 0$ for all $i$ and all $k$. We have $C_Q \cong 0$. $g_Q$ is an isomorphism and $F$ is full and faithful. \qed

We have the subsequent relevant corollary.

\begin{cor}
 Let $X$ be a quasi-projective scheme over $k$ with a perfect ample sequence $\{P_i\}$ and $F: D_{\perf}(X) \ra D_{\perf}(Y)$ an exact functor possessing left and right adjoints. If the maps
\begin{displaymath}
 \Hom(P_i,P_j[k]) \overset{\sim}{\longrightarrow} \Hom(F(P_i),F(P_j)[k])
\end{displaymath}
are isomorphisms for $i < j$ and for all $k$. Then $F$ is full and faithful on $D_{\perf}(X)$.
\end{cor}

The following is a very useful result due to Orlov, see \cite{Orl97}. We shall not recall the proof as we require no modification.

\begin{prop}
 Let $\mathcal{A}$ be an abelian category possessing an ample sequence and let $F: D^b(\mathcal{A}) \ra D^b(\mathcal{A})$ be an autoequivalence. Suppose there exists an isomorphism $f: \Id_{\Omega} \cong F|_{\Omega}$. Then $f$ can be extended to an isomorphism $\Id_{D^b(\mathcal{A})} \overset{\sim}{\longrightarrow} F$ on all of $D^b(\mathcal{A})$.
\label{lem:Orlovextendnat}
\end{prop}

\section{Resolutions of the diagonal}
\label{sec:resolvediag}
 
Let $X$ and $Y$ be projective schemes over a field $k$. Let $\mathcal{O}_X(1)$ and $\mathcal{O}_Y(1)$ denote choices of very ample sheaves on $X$ and $Y$, respectively. Then, $\mathcal{O}_X(1) \boxtimes \mathcal{O}_Y(1)$ is a very ample sheaf for $X \times Y$ since we can use $\mathcal{O}_X(1)$ and $\mathcal{O}_Y(1)$ to embed $X \times Y$ into $\P^N_k \times \P^M_k$ and then apply the Segre embedding, $\P^N_k \times \P^M_k \ra \P^{MN+M+N}_k$. The pullback of the twisting sheaf of $\P^{MN+M+N}_k$ to $X \times Y$ is $\mathcal{O}_X(1) \boxtimes \mathcal{O}_Y(1)$.

\begin{lem}
 Given any coherent sheaf $C$ on $X \times Y$, there exists a $j$ and $m$ and a surjection
\begin{displaymath}
 \left(\mathcal{O}_X(-j)\boxtimes\mathcal{O}_Y(-j)\right)^{\oplus m} \twoheadrightarrow C.
\end{displaymath}
\end{lem}

We can rewrite $\left(\mathcal{O}_X(-j)\boxtimes\mathcal{O}_Y(-j)\right)^{\oplus m}$ as $\mathcal{O}_X(-j)\boxtimes\left(\mathcal{O}_Y(-j)^{\oplus m}\right)$ and immediately deduce the following lemma.

\begin{lem}
 Any coherent sheaf $C$ on $X \times Y$ has a bounded above resolution
\begin{displaymath}
 \cdots \ra E_N \boxtimes F_N \ra \cdots \ra E_1 \boxtimes F_1 \ra E_0 \boxtimes F_0 \ra C \ra 0
\end{displaymath}
where $E_i$ are invertible sheaves on $X$ and $F_i$ are locally-free coherent sheaves on $Y$.
\end{lem}

\begin{cor}
 On $X \times X$, there is a resolution of the structure sheaf of the diagonal $\Delta X$
\begin{displaymath}
 \cdots \ra A_N \boxtimes B_N \ra \cdots \ra A_1 \boxtimes B_1 \ra A_0 \boxtimes B_0 \ra \mathcal{O}_{\Delta X} \ra 0
\end{displaymath}
where $A_i$ and $B_i$ are locally-free coherent sheaves.
\label{cor:resdiagprojective}
\end{cor}

\begin{cor}
 Assume that $X$ is quasi-projective. Then, there is a resolution of the structure sheaf of the diagonal $\Delta X$
\begin{displaymath}
 \cdots \ra A_N \boxtimes B_N \ra \cdots \ra A_1 \boxtimes B_1 \ra A_0 \boxtimes B_0 \ra \mathcal{O}_{\Delta X} \ra 0
\end{displaymath}
where $A_i$ and $B_i$ are locally-free coherent sheaves.
\label{cor:resdiagquasiproj}
\end{cor}

\proof Let $\bar{X}$ be a choice of projective closure of $X$. By corollary \ref{cor:resdiagprojective}, we have a resolution
\begin{displaymath}
 \cdots \ra A_N \boxtimes B_N \ra \cdots \ra A_1 \boxtimes B_1 \ra A_0 \boxtimes B_0 \ra \mathcal{O}_{\Delta \bar{X}} \ra 0
\end{displaymath}
of the structure sheaf of the diagonal $\Delta \bar{X}$. $A_i$ and $B_i$ are locally-free coherent sheaves on $\bar{X}$. Since $X$ is an open subset of $\bar{X}$, the restriction to $X$ is exact. So
\begin{displaymath}
 \cdots \ra A_N|_X \boxtimes B_N|_X \ra \cdots \ra A_1|_X \boxtimes B_1|_X \ra A_0|_X \boxtimes B_0|_X \ra \mathcal{O}_{\Delta X} \ra 0
\end{displaymath}
is a resolution of the diagonal in $X \times X$ with $A_i|_X$ and $B_i|_X$ locally-free. \qed

We now record and prove a useful lemma found in \cite{Kaw04}.

\begin{lem}
 Let $f:X \ra Y$ be a quasi-projective morphism with $L$ the corresponding ample sheaf. Assume that large tensor powers of $L$ have trivial higher cohomology. If $D$ is a bounded above complex of coherent sheaves on $X$ with $H^{m_0}(D)$ nonzero, then there exists an integer $N$ so that, for all $k > N$, $H^{m_0}(f_*(D \otimes L^k))$ is nonzero.
\label{lem:ssdetect}
\end{lem}

\proof Take $N$ large enough so that, for $k > N$, $H^p (f_*(H^q(D) \otimes L^k)) = 0$ for $p>0$ and $q>m_0-\dim X$ and $H^0 (f_*(H^{m_0}(D) \otimes L^k)) \not = 0$. We have a spectral sequence
\begin{displaymath}
 E_2^{pq} = H^p (f_*(H^q(D) \otimes L^k)) \Rightarrow H^{p+q}(f_*(D \otimes L^k)) 
\end{displaymath}
and we see that, thanks to our choices, $H^{m_0}(f_*(D \otimes L^k)) \not = 0$. \qed

There is one concrete and important resolution of the diagonal that we should discuss further: the Beilinson resolution of the diagonal for $\P^n_k$, \cite{Bei79}. For the polynomial algebra, $S=k[x_0,\cdots,x_n]$ we have a map of graded modules $S(-1)^{\oplus n+1} \ra S$ given by sending a basis vector $e_i$ to $x_i$. Let $M$ be the kernel of this map. The corresponding maps of coherent sheaves on $\P^n_k$,
\begin{center}
 \leavevmode
 \begin{xy}
  (-35,0)*+{0}="a"; (-23,0)*+{\tilde{M}}="b"; (0,0)*+{\mathcal{O}_{\P^n_k}(-1)^{\oplus n+1}}="c"; (22,0)*+{\mathcal{O}_{\P^n_k}}="d"; (35,0)*+{0}="e"; {\ar@{->} "a";"b"};  {\ar@{->} "b";"c"}; {\ar@{->} "c";"d"}; {\ar@{->} "d";"e"};
 \end{xy},
\end{center}
is exact. By localising to the standard open subsets, one can check that $\tilde{M} \cong \Omega_{\P^n_k/k}$. (On the affine subset where $x_i$ is nonzero, $d(x_j/x_i)$ gets mapped to $1/x_i(e_j - x_j/x_i e_i)$, and these maps glue). If we take duals and twist by $-1$, we get
\begin{center}
 \leavevmode
 \begin{xy}
  (-35,0)*+{0}="a"; (-20,0)*+{\mathcal{O}_{\P^n_k}(-1)}="b"; (0,0)*+{\mathcal{O}_{\P^n_k}^{\oplus n+1}}="c"; (20,0)*+{\mathcal{T}_{\P^n_k}(-1)}="d"; (35,0)*+{0}="e"; {\ar@{->} "a";"b"};  {\ar@{->} "b";"c"}; {\ar@{->} "c";"d"}; {\ar@{->} "d";"e"};
 \end{xy}
\end{center}
which induces an isomorphism $H^0(\P^n_k,\mathcal{O}_{\P^n_k}^{\oplus n+1}) \ra H^0(X,\mathcal{T}_{\P^n_k}(-1))$. This shows that the vector fields $\displaystyle{\frac{\del}{\del x^i}}$ are a basis of the global sections of $\mathcal{T}_{\P^n_k}$. Define a global section,
\begin{displaymath}
 s: \mathcal{O}_{\P^n_k \times \P^n_k} \ra \mathcal{O}_{\P^n_k}(1) \boxtimes \mathcal{T}_{\P^n_k}(-1),
\end{displaymath}
by setting
\begin{displaymath}
 s = \sum x_i \boxtimes \frac{\del}{\del y_i}
\end{displaymath}
where $x_i$ are coordinates on the first factor and $y_i$ are coordinates on the second. One can then check, by localizing to the affine charts, that the divisor corresponding to $s$ is exactly the diagonal. Taking the Koszul resolution associated to the section we get the resolution
\begin{gather*}
 0 \ra \mathcal{O}_{\P^n_k}(-n)\boxtimes \Omega^n_{\P^n_k}(n) \ra \mathcal{O}_{\P^n_k}(-n+1)\boxtimes \Omega^{n-1}_{\P^n_k}(n-1) \ra \cdots \\ \cdots \ra \mathcal{O}_{\P^n_k}(-1)\boxtimes \Omega_{\P^n_k}(1) \ra \mathcal{O}_{\P^n_k}\boxtimes \mathcal{O}_{\P^n_k} \ra \mathcal{O}_{\Delta} \ra 0.
\end{gather*}
Note that we can shift the degrees a bit
\begin{gather*}
 0 \ra \mathcal{O}_{\P^n_k}(m-n)\boxtimes \Omega^n_{\P^n_k}(n-m) \ra \mathcal{O}_{\P^n_k}(m-n+1)\boxtimes \Omega^{n-m-1}_{\P^n_k}(n-1) \ra \cdots \\ \cdots \ra \mathcal{O}_{\P^n_k}(m-1)\boxtimes \Omega_{\P^n_k}(1-m) \ra \mathcal{O}_{\P^n_k}(m)\boxtimes \mathcal{O}_{\P^n_k}(-m) \ra \mathcal{O}_{\Delta} \ra 0
\end{gather*}
for any $m \in \Z$.

\section{A derived Morita theorem for some quasi-projective schemes}
\label{sec:derivedMorita}
 
In this section, we assume that $X$ and $Y$ are quasi-projective over a field $k$. We also assume that $X$ possesses a line bundle $L$ which is ample and satisfies the following condition: there exists an $N$ so that, for $l > N$, $H^i(X,L^l)=0$ for $i>0$. We denote $L$ by $\mathcal{O}_X(1)$.
Any scheme that is projective over a finitely generated $k$-algebra satisfies this condition. Any scheme that is affine over a projective scheme satisfies this condition.

\begin{eg}
 One might speculate that we just need to pullback the twisting sheaf from $\P^n_k$, but this does not always work. Consider $\mathbf{A}^2_k - \{0,0\}$. Its structure sheaf is very ample, but it has higher cohomology. Since its Picard group is $\Z$, we see that, on $\mathbf{A}^2_k - \{0,0\}$, there is no ample sheaf whose large tensor powers possess trivial higher cohomology.
\end{eg}

Also, in this section, $F: D_{\perf}(X) \ra D_{\perf}(Y)$ is a full and faithful exact functor. If $X$ is not projective, we assume that $F$ admits an extension, $F: D(X) \ra D(Y)$, which commutes with coproducts.

\begin{lem}
 If $F$ has a left adjoint, then there exists $m_1$ and $m_2$ so that $F(E)$ is concentrated in $[m_1,m_2]$ for any finite rank locally-free sheaf, $E$.
\label{lem:boundedness}
\end{lem}

\proof Let us denote the left adjoint by $F^*$. Choose some embedding $Y \hookrightarrow \P^n_k$. On $\P^n_k \times \P^n_k$, we can resolve the pushforward of $\mathcal{O}_{\P^n_k}(j)$ 
\begin{gather*}
 0 \ra \mathcal{O}_{\P^n_k}(-n)\boxtimes \Omega^n_{\P^n_k}(j+n) \ra \mathcal{O}_{\P^n_k}(-n+1)\boxtimes \Omega^{n-1}_{\P^n_k}(j+n-1) \ra \cdots \\ \cdots \ra \mathcal{O}_{\P^n_k}(-1)\boxtimes \Omega_{\P^n_k}(j+1) \ra \mathcal{O}_{\P^n_k}\boxtimes \mathcal{O}_{\P^n_k}(j) \ra \mathcal{O}_{\Delta}(j) \ra 0
\end{gather*}
Applying $p_{1*}$, we see that $\mathcal{O}_{\P^n_k}(j)$ can be obtained from $\{\mathcal{O}_{\P^n_k}(-n),\cdots,\mathcal{O}_{\P^n_k}\}$ by taking shifts $[-j]$ with $0 \leq j \leq n$, direct sums, and then $n$ cones of objects within the subcategory formed by those direct sums of shifts of $\{\mathcal{O}_{\P^n_k}(-n),\ldots,\mathcal{O}_{\P^n_k}\}$. This statement pulls back to $Y$.

%

%

We choose $k_1$ and $k_2$ so that $\mathcal{H}om(F^*(\mathcal{O}_Y(j)),\mathcal{O}_X)$ is quasi-isomorphic to a complex of locally-free coherent sheaves which is zero outside $[k_1,k_2]$ for all $0 \leq j \leq n$. Take any locally-free coherent sheaf, $E$, on $Y$. Via adjunction,
\begin{gather*}
 [\mathcal{O}_Y(j),F(E)[k]] \cong [F^*(\mathcal{O}_Y(j)),E[k]] \cong [\mathcal{O}_X,\mathcal{H}om(F^*(\mathcal{O}_Y(j)),\mathcal{O}_X) \otimes E[k]]
\end{gather*}
$\mathcal{H}om(F^*(\mathcal{O}_Y(j)),\mathcal{O}_X) \otimes E$ is concentrated in $[k_1,k_2]$. Thus, $[\mathcal{O}_Y(j),F(E)[k]]$ is zero for $k$ outside $[k_1,k_2+\dim X]$. If $i$ is outside $[0,N]$, $F^*(\mathcal{O}_X(i))$ is obtained from $F^*(\mathcal{O}_X(j))$, $0 \leq j \leq N$, using a uniformly bounded number cones and uniformly bounded shifts. $[\mathcal{O}_X,\mathcal{H}om(F^*(\mathcal{O}_X(j)),\mathcal{O}_X)\otimes E [k]]$ can be computed using long exact sequences coming from the triangles needed to build $F^*(\mathcal{O}_X(j))$. Thus, we get a uniform bound, in $j$ and $k$, on $[\mathcal{O}_X,\mathcal{H}om(F^*(\mathcal{O}_Y(j)),\mathcal{O}_X) \otimes E[k]] \cong [\mathcal{O}_Y(j),F(E)[k]]$. These bounds provide our $m_1$ and $m_2$. \qed


\begin{defn}
 If the cohomologies of the image, under $F$, of complexes (in $D_{\perf}(X)$ or $D^b_{\coh}(X)$) concentrated in degree zero are uniformly bounded, we say that $F$ is \textbf{bounded}.
\end{defn}

\begin{rmk}
 Note that, if $X$ is projective, $F$ has a left adjoint if and only if $F$ extends to a functor $\tilde{F}:D^b_{\coh}(X) \ra D^b_{\coh}(Y)$. If $F$ the restriction of a functor $\tilde{F}: D^b_{\coh}(X) \ra D^b_{\coh}(Y)$, then this argument actually shows that $\tilde{F}$ is bounded.
\end{rmk}

We shall shift $F$, if necessary, and assume that the cohomology of $F(A)$, for any locally-free coherent sheaf, $A$, is concentrated in $[k_0,0]$.

Consider the following complex over $D(X \times Y)$:
\begin{displaymath}
 A_m \boxtimes F(B_m) \ra \cdots \ra A_1 \boxtimes F(B_1) \ra A_0 \boxtimes F(B_0).
\end{displaymath}
Here $A_i, B_i$ come from our choice of the resolution of the diagonal on $X \times X$.

\begin{lem}
 For any objects $C_1, C_2$ in $D_{\perf}(X)$ and $D_1, D_2$ in $D_{\perf}(Y)$,
\begin{displaymath}
 [C_1 \boxtimes D_1, C_2 \boxtimes D_2] \cong \bigoplus_{l \in \Z} [C_1,C_2[-l]] \otimes [D_1,D_2[l]].
\end{displaymath}
\end{lem}

\proof We manipulate some adjunctions:
\begin{gather*}
 [C_1 \boxtimes D_1, C_2 \boxtimes D_2] \cong [p_1^*C_1, p_1^*C_2 \otimes p_2^*(\mathcal{H}om(D_1,\mathcal{O}_Y) \otimes D_2)] \cong \\ [C_1,C_2 \otimes p_{1*}p_2^*(\mathcal{H}om(D_1,\mathcal{O}_Y) \otimes D_2))] \cong \bigoplus_{l \in \Z} [C_1,C_2[-l]] \otimes [D_1,D_2[l]].
\end{gather*}
\qed

\begin{rmk}
 As long as there are no morphisms between $D_1$ and $D_2[l]$ for $|l|$ large, we do not need perfection.
\end{rmk}

For $r>0$ and $p < q$, we have
\begin{gather*}
 [A_p \boxtimes F(B_p)[r],A_q \boxtimes F(B_q)] \cong \bigoplus_{r_1+r_2=r} [A_p[r_1],A_q] \otimes [F(B_p)[r_2],F(B_q)] \cong \\ \bigoplus_{r_1+r_2=r} [A_p[r_1],A_q] \otimes [B_p[r_2],B_q] \cong 0.
\end{gather*}
Thus a right convolution of this complex exists. Denote it by $E_m'$.

\begin{lem}
 Assume that $F(G)$ is concentrated in $[k_0,0]$ for any locally-free coherent sheaf, $G$. $H^p(E_m') = 0$ unless $p \in [m+k_0,m] \cup [k_0,0]$.
\end{lem}

\proof Let us truncate our resolution of the diagonal
\begin{displaymath}
 0 \ra T_{m-1} \ra A_m \boxtimes B_m \ra \cdots \ra A_0 \boxtimes B_0 \ra \mathcal{O}_{\Delta} \ra 0
\end{displaymath}
From our assumption on $X$, for any fixed bounded complex $D$ of coherent sheaves, there is an $N$ so that tensoring $D$ by $p_1^*\mathcal{O}_X(l)$, for $k \geq N$, and pushing forward by $p_2$ yields an exact sequence. Thus, for $l$ large, we have an exact sequence
\begin{displaymath}
 0 \ra S_{l,m-1} \ra H^0(X,A_m(l)) \boxtimes B_m \ra \cdots \ra H^0(X,A_0(l)) \boxtimes B_0 \ra \mathcal{O}_X(l) \ra 0
\end{displaymath}
where $S_{l,m-1} = p_{2*}(T_{m-1}\otimes p_1^*\mathcal{O}_X(l))$. This represents a map  $$\mathcal{O}_X(l) \ra S_{l,m-1}[m]]$$ which must be zero if we choose $m > \dim X$. Thus, the complex,
\begin{displaymath}
 H^0(X,A_m(l)) \boxtimes B_m \ra \cdots \ra H^0(X,A_0(l)) \boxtimes B_0,
\end{displaymath}
is quasi-isomorphic to $S_{l,m-1}[m] \oplus \mathcal{O}_X(l)$. The complex,
\begin{displaymath}
 H^0(X,A_m(l)) \boxtimes B_m \ra \cdots \ra H^0(X,A_0(l)) \boxtimes B_0,
\end{displaymath}
viewed as lying over $D(X)$, is convolvable. The convolution is just the complex itself as an object of $D(X)$. Hence, the convolution is quasi-isomorphic to $S_{l,m-1}[m] \oplus \mathcal{O}_X(l)$. Apply $F$ to
\begin{displaymath}
 H^0(X,A_m(l)) \boxtimes B_m \ra \cdots \ra H^0(X,A_0(l)) \boxtimes B_0
\end{displaymath}
and consider it as a complex over $D(Y)$. It is convolvable as $F$ is fully-faithful. Since exact functors map convolutions to convolutions and convolutions are unique by corollary \ref{cor:rightconvolve}, the convolution of
\begin{displaymath}
 H^0(X,A_m(l)) \otimes F(B_m) \ra  \cdots \ra H^0(X,A_0(l)) \otimes F(B_0)
\end{displaymath}
is quasi-isomorphic to $F(S_{l,m-1}[m] \oplus \mathcal{O}_X(l))$. Now, note that the complex
\begin{displaymath}
 H^0(X,A_m(l)) \otimes F(B_m) \ra  \cdots \ra H^0(X,A_0(l)) \otimes F(B_0)
\end{displaymath} also results from applying $p_{2*}(- \otimes p_1^*\mathcal{O}_X(l))$ to
\begin{displaymath}
 A_m(l) \boxtimes F(B_m) \ra \cdots \ra A_0(l) \boxtimes F(B_0).
\end{displaymath}
So
\begin{displaymath}
 \Phi_{E_m'}(\mathcal{O}_X(l)) \cong F(S_{k,m-1}[m] \oplus \mathcal{O}_X(l))
\end{displaymath}
for $l$ large(r than a fixed constant depending on $\dim X$ and $A_i,B_i$ for $0 \leq i \leq m$). Thus, the cohomology of $\Phi_{E_m'}(\mathcal{O}_X(l))$ is concentrated in $[k_0+m,m] \cup [k_0,0]$. Lemma \ref{lem:ssdetect} says that $E_m'$ has cohomology concentrated in $[k_0+m,m] \cup [k_0,0]$. \qed

We record a corollary of the proof.

\begin{cor}
 There is an $M$ (depending on $m$) so that, for any $l>M$,
\begin{displaymath}
 \Phi_{E_m'}(\mathcal{O}_X(l)) \cong F(S_{k,m-1}[m] \oplus \mathcal{O}_X(l)).
\end{displaymath}
\end{cor}

Let $\tau_{\geq k_0 -1}E_m'$ be the gentle truncation of $E_m'$. Fix $m<k_0-1$. Denote the truncation by $E$.

\begin{lem}
 There is an $N$ so that, for $l > N$, there are isomorphisms
\begin{displaymath}
 \tau_{\geq k_0 -1}(\Phi_{E_m'}(\mathcal{O}_X(l))) \cong \Phi_{E}(\mathcal{O}_X(l)).
\end{displaymath}
\end{lem}

\proof Choose $N$ large enough so that $p_{2*}(-\otimes p_1^*\mathcal{O}_X(l))$ has no higher cohomology on any element of $E_m'$ for $l > N$. Then, $p_{2*}(-\otimes p_1^*\mathcal{O}_X(l))$ commutes with truncation of $E_m'$. So $\tau_{\geq k_0-1}(\Phi_{E_m'}(\mathcal{O}_X(l))) \cong \Phi_{\tau_{\geq k_0-1} E_m'}(\mathcal{O}_X(l)) = \Phi_{E}(\mathcal{O}_X(l))$. \qed

\begin{cor}
 There exists an $N$ so that, for $l>N$, $F(\mathcal{O}_X(l))$ is quasi-isomorphic to $\Phi_E(\mathcal{O}_X(l))$. Therefore, $\Phi_E(\mathcal{O}_X(l))$ is perfect for all $l > N$.
\end{cor}

\begin{lem}
 There exists a natural isomorphism of $\Phi_E$ and $F$ on the full subcategory of $D_{\perf}(X)$ whose objects are $\{\mathcal{O}_X(l)\}_{l > N}$
\label{lem:nattransl>N}
\end{lem}

\proof We have just seen that the for each $l > N$ we have the isomorphisms of objects. The need is to make them functorial. Consider projection $\epsilon_l: F(S_{l,m-1}[m] \oplus \mathcal{O}_X(l)) \ra F(\mathcal{O}_X(l))$. Since
\begin{displaymath}
 \Hom(H^0(X,A_p(l)) \otimes F(B_p)[r],F(\mathcal{O}_X(l))) \cong \Hom(H^0(X,A_p(l)) \otimes B_p[r],\mathcal{O}_X(l)) \cong 0
\end{displaymath}
for $r > 0$ and any $p$ (both are sheaves), $\epsilon_l$ is the only morphism which makes
\begin{center}
 \leavevmode
 \begin{xy}
  (-40,10)*+{H^0(X,A_0(l)) \otimes F(B_0)}="a"; (40,10)*+{F(S_{l,m-1}[m]\oplus \mathcal{O}_X(l))}="b"; (-40,-10)*+{H^0(X,A_0(l)) \otimes F(B_0)}="c"; (40,-10)*+{F(\mathcal{O}_X(l))}="d"; {\ar@{->}^{F(r)} "a";"b"}; {\ar@{->}^= "a";"c"}; {\ar@{->}_{\epsilon_l} "b";"d"};
  {\ar@{->}^{\epsilon_lF(r)} "c";"d"};
 \end{xy}
\end{center}
commute by lemma \ref{lem:rightconvolvemorphisms}. Here $r$ is the map coming from the convolution of
\begin{displaymath}
 H^0(X,A_m(l)) \otimes B_m \ra \cdots \ra H^0(X,A_0(l)) \otimes B_0
\end{displaymath}
Similarly, the projection $\epsilon_l': \Phi_{E_m'}(\mathcal{O}_X(l) \ra \Phi_E(\mathcal{O}_X(l))$ is the only morphism making
\begin{center}
 \leavevmode
 \begin{xy}
  (-40,10)*+{H^0(X,A_0(l)) \otimes F(B_0)}="a"; (40,10)*+{\Phi_{E_m'}(\mathcal{O}_X(l))}="b"; (-40,-10)*+{H^(X,A_0(l)) \otimes F(B_0)}="c"; (40,-10)*+{\Phi_E(\mathcal{O}_X(l))}="d"; {\ar@{->}^{\Phi_{r'}(\mathcal{O}_X(l))} "a";"b"}; {\ar@{->}^= "a";"c"}; {\ar@{->}_{\epsilon_l'} "b";"d"};
  {\ar@{->}^{\epsilon_l'\Phi_{r'}(\mathcal{O}_X(l))} "c";"d"};
 \end{xy}
\end{center}
commute. Here $r'$ is the map coming from the convolution of
\begin{displaymath}
 A_m(l) \boxtimes F(B_m) \ra \cdots \ra A_0(l) \boxtimes F(B_0).
\end{displaymath}
From above, we have isomorphisms $f_l: F(\mathcal{O}_X(l)) \ra \Phi_E(\mathcal{O}_X(l))$ making 
\begin{center}
 \leavevmode
 \begin{xy}
  (-40,10)*+{H^0(X,A_0(l)) \otimes F(B_0)}="a"; (40,10)*+{F(\mathcal{O}_X(l))}="b"; (-40,-10)*+{H^0(X,A_0(l)) \otimes F(B_0)}="c"; (40,-10)*+{\Phi_E(\mathcal{O}_X(l))}="d"; {\ar@{->}^{\epsilon_lF(r)} "a";"b"}; {\ar@{->}^= "a";"c"}; {\ar@{->}_{f_l} "b";"d"}; {\ar@{->}^{\epsilon_l'\Phi_{r'}(\mathcal{O}_X(l)} "c";"d"};
 \end{xy}
\end{center}
commute. It follows from the previous observations that any map fitting in the slot of $f_l$ must be unique. If there were two such maps $f_1$ and $f_2$, then $(f_1-f_2)\epsilon_l = 0$ by uniqueness. But, $\epsilon_l$ is just a projection and we can precompose with the splitting map $F(\mathcal{O}_X(l)) \ra F(S_{l,m-1}[m] \oplus \mathcal{O}_X(l))$ to conclude that $f_1-f_2=0$. 

Now let $\alpha: \mathcal{O}_X(l_1) \ra \mathcal{O}_X(l_2)$ be a morphism for $l_1,l_2 > N$. Then there is uniquely determined $g$ so that
\begin{center}
 \leavevmode
 \begin{xy}
  (-40,10)*+{H^0(X,A_0(l_1)) \otimes F(B_0)}="a"; (40,-10)*+{\Phi_{E}(\mathcal{O}_X(l_2))}="b"; (-40,-10)*+{H^0(X,A_0(l_2)) \otimes F(B_0)}="c"; (40,10)*+{F(S_{l_1,m-1}[m]\oplus \mathcal{O}_X(l_1))}="d"; {\ar@{->}^{\epsilon_{l_1}'\Phi_{r'}(\mathcal{O}_X(l_2))} "c";"b"}; {\ar@{->}^{H^0(X,\alpha) \otimes \Id} "a";"c"}; {\ar@{->}_{g} "d";"b"}; {\ar@{->}^{F(r)} "a";"d"};
 \end{xy}
\end{center}
commutes. Both $\Phi_E(\alpha)f_{l_1}\epsilon_{l_1}$ and $f_{l_2}F(\alpha) \epsilon_{l_1}$ fit. Therefore
\begin{displaymath}
 \Phi_E(\alpha)f_{l_1}\epsilon_{l_1} = f_{l_2}F(\alpha) \epsilon_{l_1} = g
\end{displaymath}
and consequently $\Phi_E(\alpha)f_{l_1} = f_{l_2}F(\alpha)$. \qed

\begin{lem}
 Let $F: D_{\perf}(X) \ra D_{\perf}(Y)$ be an equivalence and let $E$ be as constructed above. Assume there is natural isomorphism of $F$ with $\Phi_E$ on the full subcategory formed by $\{\mathcal{O}_X(l)\}_{l > N}$. Then, $\Phi_E: D(X) \ra D(Y)$ is an equivalence.
\end{lem}

\proof Let $\mathcal{S}_L$ be the smallest full subcategory of $D(X)$ consisting of objects $B$ for which $[\mathcal{O}_X(l),B] \ra [\Phi_E(\mathcal{O}_X(l)),\Phi_E(B)]$ is a bijection for $l > N$. $\mathcal{S}_L$ is triangulated, closed under coproducts (as $\Phi_E$ takes perfect objects to perfect objects), and contains $\{\mathcal{O}_X(l)\}_{l > N}$. Thus, $\mathcal{S}_L \cong D(X)$. Let $\mathcal{S}_R$ be the smallest full subcategory of $D(X)$ consisting of objects $A$ for which $[A,B] \ra [\Phi_E(A),\Phi_E(B)]$ is a bijection for all $B \in D(X)$. $\mathcal{S}_R$ is triangulated, closed under coproducts (naturally), and contains $\{\mathcal{O}_X(l)\}_{l > N}$. Thus, $\mathcal{S}_R \cong D(X)$. Consequently, $\Phi_E$ is full and faithful. Since $F$ is an equivalence, the smallest triangulated subcategory containing $\{F(\mathcal{O}_X(l))\}_{l > N}$ is all of $D_{\perf}(Y)$. Therefore, the essential image of $\Phi_E$ contains all of $D_{\perf}(Y)$. The essential image is closed under triangles and coproducts and is therefore all of $D(Y)$. \qed

We can now state a derived Morita theorem (in the sense of Rickard) for certain quasi-projective schemes. For the original statement, see \cite{Ric89}.

\begin{thm}
 Let $X$ and $Y$ be quasi-projective schemes over a field $k$. Assume that $X$ possesses an ample line bundle, sufficiently high powers of which have trivial higher cohomology. The following are equivalent:
\begin{itemize}
 \item There is an exact equivalence $D(X) \ra D(Y)$.
 \item There is an exact equivalence $D^b(X) \ra D^b(Y)$.
 \item There is an object $E \in D^b_{\coh}(X \times Y)$ so that $\Phi_E : D(X) \ra D(Y)$ is an exact equivalence.
\end{itemize}
 Moreover, if $X$ is projective, then the following are also equivalent:
\begin{itemize}
 \item There is an exact equivalence $D_{\perf}(X) \ra D_{\perf}(Y)$.
 \item There is an exact equivalence $D^b_{\coh}(X) \ra D^b_{\coh}(Y)$.
\end{itemize}
\end{thm}

\proof In the quasi-projective case, the work of this section showed the equivalence of the first and third conditions. Tracing out the arguments, we see that we can restrict ourselves to bounded complexes in the arguments and still conclude that the third condition holds.

In the projective case, we can restrict our attention to perfect objects to conclude that the old third condition holds. The fourth and fifth conditions are equivalent by locally-finite duality. \qed

Orlov's original theorem says something a bit stronger. It says that $F$ and $\Phi_E$ are isomorphic on all of $D^b_{\coh}(X)$ if $X$ and $Y$ are smooth and projective. In the next section, we push a little harder and provide an extension of Orlov's result.

\section{Equivalences and Fourier-Mukai transforms}
\label{sec:repfunct}
 
As in last section, $F$ is a full and faithful exact functor $D_{\perf}(X) \ra D_{\perf}(Y)$. However, in this section, we take $X$ and $Y$ to be projective schemes over a field.

We are almost ready to deduce an extension of Orlov's result. We first need to know that $\Phi_E$ takes bounded complexes of coherent sheaves to bounded complexes of coherent sheaves.

\begin{lem}
 $E$ is independent of the choice of $m$ if $m < k_0 -2$.
\label{lem:EishocolimE_m}
\end{lem}

\proof Instead of fixing $m$, we shall now let $m$ run starting from $m < k_0-2$. From lemma \ref{lem:rightconvolvemorphisms}, there is a morphism $\phi_m: E_m' \ra E_{m-1}'$. $\phi_m$ induces maps $\Phi_{\phi_m}(\mathcal{O}_X(l)): \Phi_{E_m'}(\mathcal{O}_X(l)) \ra \Phi_{E_{m-1}'}(\mathcal{O}_X(l))$. Applying lemma \ref{lem:rightconvolvemorphisms} again, we see that this corresponds to a morphism $F(S_{l,m-1}[-m] \oplus \mathcal{O}_X(l)) \ra F(S_{l,m-2}[-m+1] \oplus \mathcal{O}_X(l))$. It is easy to check that the morphism, $F(S_{l,m-1}[-m] \oplus \mathcal{O}_X(l)) \ra F(S_{l,m-2}[-m+1] \oplus \mathcal{O}_X(l))$, comes from applying $F$ to the morphism $\lambda_m \oplus \id_{\mathcal{O}_X(l)}: S_{l,m-1}[-m] \oplus \mathcal{O}_X(l) \ra S_{l,m-2}[-m-1] \oplus \mathcal{O}_X(l)$ where $\lambda_m: S_{l,m-1}[-m] \ra S_{l,m-2}[-m+1]$ is the extension corresponding to the map, $H^0(X,A_{m-1}(l)) \otimes F(B_{m-1}) \ra H^0(X,A_m(l)) \otimes F(B_m)$. Since this is true for all $l \gg 0$, we see that $\Phi_{\phi_m}(\mathcal{O}_X(l))$ is a quasi-isomorphism in degrees $>m+1$. We truncate $E_m'$ above $k_0-1$. Denote the result by $E_m$. We have a diagram
\begin{center}
 \leavevmode
 \begin{xy}
  (-45,10)*+{\tau_{\leq k_0-1}E_m'}="a"; (-15,10)*+{E_m'}="b"; (15,10)*+{\tau_{\geq k_0-1}E_m'}="c"; (45,10)*+{\tau_{\leq k_0-1}E_m'[1]}="d"; {\ar@{->} "a";"b"}; {\ar@{->} "b";"c"}; {\ar@{->} "c";"d"};
  (-45,-10)*+{\tau_{\leq k_0-1}E_{m-1}'}="a'"; (-15,-10)*+{E_{m-1}'}="b'"; (15,-10)*+{\tau_{\geq k_0-1}E_{m-1}'}="c'"; (45,-10)*+{\tau_{\leq k_0-1}E_{m-1}'[1]}="d'"; {\ar@{->} "a'";"b'"}; {\ar@{->} "b'";"c'"}; {\ar@{->} "c'";"d'"};
  {\ar@{->}^{\phi_m} "b";"b'"};
 \end{xy}
\end{center}
which can be completed to a morphism of triangles
\begin{center}
 \leavevmode
 \begin{xy}
  (-45,10)*+{\tau_{\leq k_0-1}E_m'}="a"; (-15,10)*+{E_m'}="b"; (15,10)*+{\tau_{\geq k_0-1}E_m'}="c"; (45,10)*+{\tau_{\leq k_0-1}E_m'[1]}="d"; {\ar@{->} "a";"b"}; {\ar@{->} "b";"c"}; {\ar@{->} "c";"d"};
  (-45,-10)*+{\tau_{\leq k_0-1}E_{m-1}'}="a'"; (-15,-10)*+{E_{m-1}'}="b'"; (15,-10)*+{\tau_{\geq k_0-1}E_{m-1}'}="c'"; (45,-10)*+{\tau_{\leq k_0-1}E_{m-1}'[1]}="d'"; {\ar@{->} "a'";"b'"}; {\ar@{->} "b'";"c'"}; {\ar@{->} "c'";"d'"};
  {\ar@{->} "a";"a'"}; {\ar@{->}^{\phi_m} "b";"b'"}; {\ar@{->}^{\psi_m} "c";"c'"}; {\ar@{->} "d";"d'"};
 \end{xy}
\end{center}
since $[\tau_{\leq k_0-1}E_m',\tau_{\geq k_0-1}E_{m-1}']=0$ ($m < k_0 - 2$). As $[\tau_{\leq k_0-1}E_m'[1],\tau_{\geq k_0-1}E_{m-1}']=0$, $\psi_m$ is unique. If we apply $\Phi_{-}(\mathcal{O}_X(l))$ to this diagram we get
\begin{center}
 \leavevmode
 \resizebox{\columnwidth}{!}{\begin{xy}
  (-60,10)*+{F(S_{l,m-1}[m])}="a"; (-20,10)*+{F(S_{l,m-1}[m]\oplus \mathcal{O}_X(l))}="b"; (20,10)*+{F(\mathcal{O}_X(l))}="c"; (60,10)*+{F(S_{l,m-1}[m+1])}="d"; {\ar@{->} "a";"b"}; {\ar@{->} "b";"c"}; {\ar@{->} "c";"d"};
  (-60,-10)*+{F(S_{l,m-2}[m-1])}="a'"; (-20,-10)*+{F(S_{l,m-1}[m]\oplus \mathcal{O}_X(l))}="b'"; (20,-10)*+{F(\mathcal{O}_X(l))}="c'"; (60,-10)*+{F(S_{l,m-2}[m])}="d'"; {\ar@{->} "a'";"b'"}; {\ar@{->} "b'";"c'"}; {\ar@{->} "c'";"d'"};
  {\ar@{->} "a";"a'"}; {\ar@{->}^{\Phi_{\phi_m}(\mathcal{O}_X(l))} "b";"b'"}; {\ar@{->}^{\Phi_{\psi_m}(\mathcal{O}_X(l))} "c";"c'"}; {\ar@{->} "d";"d'"};
 \end{xy}}
\end{center}
The map fitting into the slot of $\Phi_{\psi_m}(\mathcal{O}_X(l))$ is unique, and the composition of the inclusion of $F(\mathcal{O}_X(l))$, $\Phi_{\phi_m}(\mathcal{O}_X(l))$, and projection onto $F(\mathcal{O}_X(l))$ also fits there. Thus, $\Phi_{\psi_m}(\mathcal{O}_X(l))$ is a quasi-isomorphism for $l \gg 0$. By lemma \ref{lem:ssdetect}, $\psi_m$ is an quasi-isomorphism for each $m$. \qed

\begin{lem}
 Assume that $F$ admits an extension, $F: D^b_{\coh}(X) \ra D^b_{\coh}(Y)$ and that $F$ is bounded (with $[k_0,0]$ being the region with nonvanishing cohomology). Then, $\Phi_{E_m'}(C)$ has cohomology concentrated in $[m+k_0,m+\dim X] \cup [k_0,\dim X]$ for any coherent sheaf $C$.
\end{lem}

\proof Consider the complex
\begin{displaymath}
 (A_m \otimes C) \boxtimes F(B_m) \ra \cdots \ra (A_0 \otimes C) \boxtimes F(B_0).
\end{displaymath}
Since $A_i,F(B_i)$ are perfect, this complex can be convolved, and the convolution is quasi-isomorphic to $E_m'\otimes p_1^*C$.
Let $F_m$ denote the convolution of the complex
\begin{displaymath}
 A_m \boxtimes B_m  \ra \cdots \ra A_0 \boxtimes B_0.
\end{displaymath}
Then, $F_m \otimes p_1^*C$ is the convolution of the complex
\begin{displaymath}
 (A_m \otimes C) \boxtimes B_m  \ra \cdots \ra (A_0 \otimes C) \boxtimes B_0
\end{displaymath}
Let $R_{m-1}$ denote the kernel of the map $(A_{m-1} \otimes C) \boxtimes B_{m-1} \ra (A_{m} \otimes C) \boxtimes B_{m}$. Choose $N$ large enough so that $A_i \otimes C(l)$ has no higher cohomology for $l > N$ and $i \geq m$ and $p_{2*}R_{l,m-1}=p_{2*}(p_1^*\mathcal{O}_X(l) \otimes R_{m-1})$ is concentrated in degree zero. We have an triangle
\begin{center}
 \leavevmode
 \begin{xy}
  (-20,10)*+{p_{2*}(R_{l,m-1}[m])}="a"; (20,10)*+{C}="b"; (0,-5)*+{H^0(X,A_m\otimes C(l))\otimes B_m \ra \cdots \ra H^0(X,A_0\otimes C(l))\otimes B_0}="c"; {\ar@{-->} "b";"a"}; {\ar@{->} "c";"b"}; {\ar@{->} "a";"c"}
 \end{xy}
\end{center}
The convolution of the complex
\begin{displaymath}
 F(H^0(X,A_m\otimes C(l))\otimes B_m) \ra \cdots \ra F(H^0(X,A_0\otimes C(l))\otimes B_0)
\end{displaymath}
is quasi-isomorphic to $\Phi_{E_m'}(C(l))$ by uniqueness of convolutions.
Applying $F$ to the previous triangle then gives the triangle
\begin{center}
 \leavevmode
 \begin{xy}
  (-20,10)*+{F(p_{2*}(R_{l,m-1}))[m]}="a"; (20,10)*+{F(C)}="b"; (0,-5)*+{\Phi_{E_m'}(C(l))}="c"; {\ar@{-->} "b";"a"}; {\ar@{->} "c";"b"}; {\ar@{->} "a";"c"}
 \end{xy}
\end{center}
Applying $\mathcal{H}om(\mathcal{O}_X,-)$ gives us a long exact sequence of cohomology sheaves. Since $F$ is bounded, $\Phi_{E_m'}(C(l))$ is concentrated in $[k_0+m,m] \cup [k_0,0]$. Applying lemma \ref{lem:ssdetect}, we get that $E_m' \otimes p_1^*C$ is concentrated in $[k_0+m,m] \cup [k_0,0]$ and the claim follows. \qed

\begin{lem}
 Assume that $F$ admits an extension, $F: D^b_{\coh}(X) \ra D^b_{\coh}(Y)$ and that $F$ is bounded. Then, $\Phi_E$ takes $D^b_{\coh}(X)$ to $D^b_{\coh}(Y)$.
\end{lem}

\proof We have the triangle
\begin{center}
 \leavevmode
 \begin{xy}
  (-20,10)*+{\tau_{\leq k_0-1}E_m'}="a"; (20,10)*+{E_m'}="b"; (0,-5)*+{\tau_{\geq k_0-1}E_m'}="c"; {\ar@{->} "a";"b"}; {\ar@{->} "b";"c"}; {\ar@{-->} "c";"a"}
 \end{xy}
\end{center}
to which we apply $\Phi_-(C)$:
\begin{center}
 \leavevmode
 \begin{xy}
  (-20,10)*+{\Phi_{\tau_{\leq k_0-1}E_m'}(C)}="a"; (20,10)*+{\Phi_{E_m'}(C)}="b"; (0,-5)*+{\Phi_{E}(C)}="c"; {\ar@{->} "a";"b"}; {\ar@{->} "b";"c"}; {\ar@{-->} "c";"a"}
 \end{xy}
\end{center}
By the previous lemma, $\Phi_{E_m'}(C)$ is concentrated in $[m+k_0,m+\dim X] \cup [k_0,\dim X]$. $\Phi_{\tau_{\leq k_0-1}E_m'}(C)$ is concentrated in $(-\infty,m+\dim X]$. $\Phi_{E}(C)$ is independent of the choice of $m$ for $m$ large. By considering the long exact sequence of cohomology sheaves as $m$ grows large, we see $\Phi_E(C)$ must be concentrated in $[k_0,\dim X]$. Consequently, $\Phi_E$ must take an bounded complex of coherent sheaves to another bounded complex of coherent sheaves. \qed

We next extend the natural isomorphism found in lemma \ref{lem:nattransl>N}.

\begin{lem}
 There exists a natural isomorphism between $\Phi_E|_{\Omega}$ and $F|_{\Omega}$ where $\Omega$ is the full subcategory formed by $\mathcal{O}_X(l)$ for $l \in \Z$.
\end{lem}

\proof We proceed by downward induction. Choose an embedding of $X$ in $\P^n_k$. Then, we have a exact sequence
\begin{displaymath}
 0 \ra \mathcal{O}_X \ra V_n \otimes \mathcal{O}_X(1) \ra \cdots \ra V_0 \otimes \mathcal{O}_X(n+1) \ra 0
\end{displaymath}
where $V_p = H^n(\P^n_k,\Omega^p_{\P^n_k}(p-n-1))$. This comes from the Beilinson resolution of the diagonal. Twisting we have an exact sequence
\begin{displaymath}
 0 \ra \mathcal{O}_X(k) \ra V_n \otimes \mathcal{O}_X(k+1) \ra \cdots \ra V_0 \otimes \mathcal{O}_X(k+n+1) \ra 0
\end{displaymath}
Note that $\mathcal{O}_X(k)$ is the left convolution of the complex
\begin{displaymath}
 V_n \otimes \mathcal{O}_X(k+1) \ra \cdots \ra V_0 \otimes \mathcal{O}_X(k+n+1)
\end{displaymath}
We already know that we have a natural isomorphism $\Phi_E \ra F$ on the subcategory formed by $\mathcal{O}_X(l)$ for $l > N$. If we set $k=N$ above and use the natural transformation, we get a morphism of complexes
\begin{center}
 \leavevmode
 \resizebox{\columnwidth}{!}{\begin{xy}
  (-55,10)*+{V_n \otimes \Phi_E(\mathcal{O}_X(N+1))}="a"; (-10,10)*+{V_{n-1} \otimes \Phi_E(\mathcal{O}_X(N+2))}="b"; (20,10)*+{\cdots}="c"; (50,10)*+{V_0 \otimes \Phi_E(\mathcal{O}_X(N+n+1))}="e"; {\ar@{->} "a";"b"}; {\ar@{->} "b";"c"}; {\ar@{->} "c";"e"};
  (-55,-10)*+{V_n \otimes F(\mathcal{O}_X(N+1))}="a'"; (-10,-10)*+{V_{n-1} \otimes F(\mathcal{O}_X(N+2))}="b'"; (20,-10)*+{\cdots}="c'"; (50,-10)*+{V_0 \otimes F(\mathcal{O}_X(N+n+1))}="e'"; {\ar@{->} "a'";"b'"}; {\ar@{->} "b'";"c'"}; {\ar@{->} "c'";"e'"};
  {\ar@{->}^{\Id_{V_n}\otimes f_{N+1}} "a";"a'"}; {\ar@{->}^{\Id_{V_{n-1}}\otimes f_{N+2}} "b";"b'"}; {\ar@{->}^{\Id_{V_0}\otimes f_{N+n+1}} "e";"e'"};
 \end{xy}}
\end{center}
which gives a unique morphism $f_{N}: \Phi_E(\mathcal{O}_X(N)) \ra F(\mathcal{O}_X(N))$ by lemma \ref{lem:leftconvolvemorphisms}. Working downward, we get all $f_l$ for $l \in \Z$. The morphism of the complexes is an isomorphism as each $f_l$ is an isomorphism. 

To check that these are natural, we take a morphism $\alpha: \mathcal{O}_X(l) \ra \mathcal{O}_X(l')$. It induces a morphism of complexes
\begin{center}
 \leavevmode
 \begin{xy}
  (-30,10)*+{V_n \otimes \Phi_E(\mathcal{O}_X(l+1))}="a"; (0,10)*+{\cdots}="b"; (30,10)*+{V_0 \otimes \Phi_E(\mathcal{O}_X(l+n+1))}="c"; {\ar@{->} "a";"b"}; {\ar@{->} "b";"c"};
  (-30,-10)*+{V_n \otimes F(\mathcal{O}_X(l'+1))}="a'"; (0,-10)*+{\cdots}="b'"; (30,-10)*+{V_0 \otimes F(\mathcal{O}_X(l'+n+1))}="c'"; {\ar@{->} "a'";"b'"}; {\ar@{->} "b'";"c'"};
  {\ar@{->}^{\Id_{V_n}\otimes f_{l'+1}\Phi_E(\alpha \otimes \id_{\mathcal{O}_X(1)})} "a";"a'"}; {\ar@{->}^{\Id_{V_0}\otimes f_{l'+n+1}\Phi_E(\alpha \otimes \id_{\mathcal{O}_X(n+1)})} "c";"c'"};
 \end{xy}
\end{center}
which corresponds to a unique morphism $\Phi_E(\mathcal{O}_X(l)) \ra F(\mathcal{O}_X(l'))$. Both $F(\alpha)f_l$ and $f_{l'}\Phi_E(\alpha)$ fit into the diagram given in lemma \ref{lem:leftconvolvemorphisms}. Thus, they are equal and the isomorphism between $\Phi_E$ and $F$ on the subcategory consisting of $\{\mathcal{O}_X(l)\}_{l \in \Z}$ is natural. \qed

We now need a useful lemma:

\begin{lem}
 Let $F: D_{\perf}(X) \ra D_{\perf}(Y)$ be a full and faithful exact functor with left and right adjoints and $G: D^b_{\coh}(X) \ra D^b_{\coh}(Y)$ another exact functor with a left pseudo-adjoint. Let $\Omega = \{ P_k \}$ is a perfect ample sequence in $D_{\perf}(X)$.  Assume $G$ preserves perfection, and there is an isomorphism of functors $f_{\Omega}: F|_{\Omega} \ra G|_{\Omega}$. Then there exists an isomorphism of functors $f: F|_{D_{\perf}(X)} \ra G|_{D_{\perf}(Y)}$ which extends $f_{\Omega}$. Moreover, if $F$ is an equivalence, then $F$ and $G$ are naturally isomorphic.
\label{lem:key}
\end{lem}

\proof Since $F$ posseses a left adjoint, the pseudo-right adjoint of $F$ is an extension of $F$ to the bounded derived categories. We will denote this extension by $F$ also. Recall that a perfect ample sequence must split generate $D_{\perf}(X)$. Thus, $G$ must take perfect objects to perfect objects. The restriction of $G$ to $D_{\perf}(X)$ has a right pseudo-adjoint which we will simply denote $G^{\vee}$.

We have a natural morphism of functors, $\Id_{D_{\perf}(X)} \ra F^{\vee}F$, which is an isomorphism as $F$ is full and faithful. Since $G$ is isomorphic to $F$ when restricted to $\Omega$, it also full and faithful on $D_{\perf}(X)$ and the natural morphism of functors, $\leftexp{\vee}{G} G|_{D_{\perf}(X)} \ra \Id_{D_{\perf}(X)}$, is an isomorphism.

$\leftexp{\vee}{G} F$ is the left adjoint to $F^{\vee} G|_{D_{\perf}(X)}$ on $D_{\perf}(X)$. They are both isomorphic to the identity when restricted to $\Omega$ so they are both fully faithful. Since they are adjoints to each other, they are quasi-inverses. $F^{\vee} G$ is right pseudo-adjoint (on $D^b_{\coh}(X)$) to $\leftexp{\vee}{G} F$. Thus, $F^{\vee} G$ is an autoequivalence of $D^b_{\coh}(X)$. By proposition \ref{lem:Orlovextendnat}, there is an natural isomorphism $\Id_{D^b_{\coh}(X)} \ra F^{\vee} G$. Restricting to $D_{\perf}(X)$ and using adjunction, we obtain a morphism of functors $f: F|_{D_{\perf}(X)} \ra G|_{D_{\perf}(X)}$.

Take any object $a$ of $D_{\perf}(X)$. Let
\begin{center}
 \leavevmode
 \begin{xy}
  (-10,10)*+{F a}="a"; (10,10)*+{G a}="b"; (0,-5)*+{c}="c"; {\ar@{->}^{f_a} "a";"b"}; {\ar@{->} "b";"c"}; {\ar@{->}^{[1]} "c";"a"}
 \end{xy}
\end{center}
be a distinguished triangle. Since $F^{\vee}(f_a)$ is an isomorphism, we have $F^{\vee}c \cong 0$. From
\begin{displaymath}
 \Hom(\omega,G^{\vee}c) \cong \Hom(G\omega,c) \cong \Hom(F\omega,c) \cong \Hom(\omega,F^{\vee}c) \cong 0
\end{displaymath}
for any $\omega \in \Omega$, we see that $G^{\vee}c \cong 0$. Hence $\Hom(G(a),c) \cong 0$, thus $Fa \cong Ga \oplus c[-1]$. But
\begin{displaymath}
 \Hom(Fa,c[-1]) \cong \Hom(a,F^{\vee}c[-1]) \cong 0
\end{displaymath}
hence $c \cong 0$. Hence, $f$ is a natural isomorphism.

Assume that $F$ is an auto-equivalence. Note that by uniqueness $\leftexp{\vee}{F} \cong F^{-1}$ and $\leftexp{\vee}{G} \cong G|_{D_{\perf}(X)}^{-1}$. We have an natural isomorphism $g: F|_{D_{\perf}(X)}^{-1} \ra G|_{D_{\perf}(X)}^{-1}$ which corresponds to $f: F|_{D_{\perf}(X)} \ra G|_{D_{\perf}(Y)}$. Taking the right dual of this natural transformation we get a natural isomorphism $g^{\vee}: \left(\leftexp{\vee}{F}\right)^{\vee} \ra \left(\leftexp{\vee}{G}\right)^{\vee}$ which gives the isomorphism between $F$ and $G$.
\qed

If $X$ and $Y$ are projective, we have seen that $\Phi_E$ must have left and right pseudo-adjoints because it preserves perfection and bounded coherence.

\begin{thm}
 Let $X$ and $Y$ be projective schemes over a field $k$. If $F: D_{\perf}(X) \ra D_{\perf}(Y)$ is a full and faithful functor with a left and a right adjoint, then $F$ is isomorphic to the restriction of $\Phi_E$ to $D_{\perf}(X)$ for $E \in D^b_{\coh}(X \times Y)$.
\end{thm}

A particular case of this is:

\begin{cor}
 Let $X$ and $Y$ be projective schemes over a field $k$. Let $F: D_{\perf}(X) \ra D_{\perf}(Y)$ be an exact equivalence. There exists an $E \in D^b_{\coh}(X \times Y)$ and a natural isomorphism $F \cong \Phi_E|_{D_{\perf}(X)}$.
\end{cor}

We can give another corollary that is a consequence of locally-finite duality.

\begin{cor}
 Let $X$ and $Y$ be projective schemes over a field $k$. If $F: D^b_{\coh}(X) \ra D^b_{\coh}(Y)$ is an exact equivalence, then there exists an $E \in D^b_{\coh}(X \times Y)$ and a natural isomorphism $F \cong \Phi_E|_{D^b_{\coh}(X)}$.
\end{cor}

\begin{rmk}
 The complex, $E$, appearing in the previous results is unique.
\end{rmk}

\proof Note that the complex, $E$, given in these results is unique. Suppose we have two complexes, $E$ and $\tilde{E}$, and natural isomorphisms between $F$ and both $\Phi_E$ and $\Phi_{\tilde{E}}$. The complexes
\begin{displaymath}
 A_m \boxtimes \Phi_E(B_m) \ra \cdots \ra A_0 \boxtimes \Phi_E(B_0)
\end{displaymath}
and
\begin{displaymath}
 A_m \boxtimes \Phi_{\tilde{E}}(B_m) \ra \cdots \ra A_0 \boxtimes \Phi_{\tilde{E}}(B_0)
\end{displaymath}
are isomorphic. We can take $A_0$ and $B_0$ to be $\mathcal{O}_X$. Using the counit $p_2^*p_{2*} \ra \Id$, we get a morphism $A_0 \boxtimes \Phi_{\tilde{E}}(B_0) \ra \tilde{E}$. Using the isomorphism of $\Phi_E$ and $\Phi_{\tilde{E}}$, we can view this as a map from the complex $A_i \boxtimes \Phi_E(B_i)$ to $\tilde{E}$. Consequently, we get a map between their convolutions $E_m' \ra \tilde{E}$. The induced map $\Phi_{E_m}(\mathcal{O}_X(l)) \ra \Phi_{\tilde{E}}(\mathcal{O}_X(l))$ is a quasi-isomorphism in degrees $> m+\dim X+\dim Y$. Consequently, the map $E_m \ra \tilde{E}$ is a quasi-isomorphism in degrees $> m+\dim X+\dim Y$. After truncation, we get a quasi-isomorphism between $\tilde{E}$ and $E$. \qed

\bibliographystyle{plain}
\bibliography{thesis}
\end{document}